\documentclass[preprint,10pt]{elsarticle}
\usepackage{amsthm,amsmath,amssymb,mathrsfs}
\usepackage[colorlinks=true,citecolor=black,linkcolor=black,urlcolor=blue]{hyperref}
\usepackage{graphicx}
\usepackage{float}
\usepackage{epstopdf}
\usepackage{epsfig}
\usepackage{extarrows,chngpage,array,float,eqnarray}     %长等号的宏包

\journal{}

\theoremstyle{plain}
\newtheorem{theorem}{Theorem}[section]

\theoremstyle{definition}
\newtheorem{definition}[theorem]{Definition}

\newtheorem{proposition}[theorem]{Proposition}
\newtheorem{question}[theorem]{Question}

\theoremstyle{remark}

\allowdisplaybreaks

\title{A recursive definition for the polymatroid Tutte polynomial}

\author{Xiaxia Guan$^{a}$,~~Xian'an Jin$^{b}$,~~Weiling Yang$^{b}$\footnote{Corresponding author.} \\
\small $^a$Department of Mathematics, Taiyuan University of Technology, P. R. China\\
\small $^b$School of Mathematical Sciences, Xiamen University, P. R. China\\
\small \emph{Email addresses}: guanxiaxia@tyut.edu.cn; xajin@xmu.edu.cn;  ywlxmu@163.com}
\date{}
\begin{document}
\begin{abstract}
The Tutte polynomial is a significant invariant of  graphs and matroids.  It is well-known that it has three equivalent definitions: bases expansion, rank generating function, and deletion-contraction formula. The polymatroid Tutte polynomial $\mathscr{T}_{P}$ generalizes the Tutte polynomial from matroids to polymatroids $P$. In \emph{[Adv. Math. 402 (2022) 108355.]} and \emph{[J. Combin. Theory Ser. A 188 (2022) 105584]}, the authors provided bases expansion and rank generating function constructions for $\mathscr{T}_{P}$, respectively. In \emph{[Int. Math. Res. Not. 19 (2025) rnaf302]}, a  recursive formula for $\mathscr{T}_{P}$ was obtained. In this paper, we show that the recursive formula itself can be used to define the polymatroid Tutte polynomial independently.
\end{abstract}


\maketitle
\section{Introduction}
\noindent

The Tutte polynomial \cite{Tutte} is a well-studied topic in graph theory and matroid theory. (Crapo \cite{Crapo} extended the Tutte polynomial from graphs to matroids.) It contains a large number of polynomials as special cases, such as chromatic polynomials and flow polynomials in graph theory, characteristic  polynomials in matroid theory, and Jones polynomial in knot theory and so on. We start with the definition of matroids.

\begin{definition}\cite{Whitney}
  Let $E$ be a finite set. Let $r:2^{E}\rightarrow \mathbb{Z}_{\geq 0}$  be a function satisfying  the following three conditions:
        \begin{enumerate}
\item[(r1)] $0\leq r(E')\leq |E'|$ for any $E'\subseteq E$;

\item[(r2)] $r(E')\leq r(E'')$ for any $E'\subseteq E''\subseteq E$;

\item[(r3)] $r(E')+r(E'')\geq r(E'\cup E'')+r(E'\cap E'')$ for any $E',E''\subseteq E$.
\end{enumerate}
 Define
$\mathcal{I}:=\{E'\subseteq E||E'|=r(E')\}$.
Then $M=(E,\mathcal{I})$ is called a \emph{matroid}, where $E$ is called the \emph{ground set} of $M$,
 the function $r$ is called the \emph{rank function} of $M$.
Any subset $I\in \mathcal{I}$ is called a \emph{basis} of $M$ if $r(I)=r(E)$.
\end{definition}
Tutte polynomials of matroids have the following three equivalent definitions.
\begin{definition} (Bases expansion)\label{Bases}
Let $M$ be a matroid over $[n]$. Order the elemants of $E$. Let $B$ be a basis of $M$.
Then  $i\in B$  is called \emph{internally active} if there is no $j<i$ such that $(B-i)\cup j$  is
a basis of $M$. Let $i_{M}(A)$  denote the number of internally active elements with respect to  $B$.
An element $i\in [n]\setminus B$ is called \emph{externally active} if there is no $j<i$ such that $(B \cup i)-j$ is
a basis of $M$. Let $j_{M}(B)$ denote the number of externally active elements with respect to  $B$.

The \emph{Tutte polynomial} $T_{M}(x, y)$ of $M$ is defined as
$$T_{M}(x, y)=\sum_{B \text{ is a basis of }M}x^{i_{M}(B)}y^{j_{M}(B)}.$$
\end{definition}

\begin{definition} (Rank generating function)\label{Rank}
Let $M$ be a matroid over $E$ with rank function $r$. Then
\begin{eqnarray*}
T_M(x,y)=\sum_{A\subseteq E}(x-1)^{r(E)-r(A)}(y-1)^{|A|-r(A)}.
\end{eqnarray*}
\end{definition}

\begin{definition} (Deletion-contraction formula) \label{Reduction}
Let $M$ be a matroid over $E$ with rank function $r$ and $e\in E$. Then  the \emph{deletion} $M\setminus e$ and \emph{contraction} $M/e$, which are matroids on $E\setminus e$, are given by the rank functions $r_{M\setminus e}(A)=r(A)$ and $r_{M/e}(A)=r(A\cup e)-r(e)$, for any subset $A\subseteq E\setminus e$, respectively. In particular, $e\in E$ is a \emph{loop} of $M$ if $r(e)=0$.
Any $e\in E$ is a \emph{coloop} of $M$ if $r(E\setminus e)=r(E)-1$.

The \emph{Tutte polynomial} $T_{M}(x, y)$ of $M$ is defined as follows. If $E=\emptyset$, then $T_M(x,y)=1$. Let $e\in E$. Then
\begin{equation*}
T_{M}(x,y)=\left\{\begin{array}{ll}
xT_{M/e}(x,y)&\text{if $e$ is a coloop};\\
yT_{M\setminus e}(x,y)&\text{if $e$ is a loop};\\
T_{M/e}(x,y)+T_{M\setminus e}(x,y)&\text{otherwise}.
\end{array}\right.
\end{equation*}
\end{definition}

The order of $E$ plays an implicit role in Definitions \ref{Bases} and \ref{Reduction}, but Definition \ref{Rank} depends only on matroids and not on this order. Hence, the Tutte polynomial depends only on matroids.

Let  $[n]=\{1,2,\ldots,n\}$, $2^{[n]}=\{I|I\subseteq[n]\}$, and let $\textbf{e}_{1},\textbf{e}_{2},\ldots,\textbf{e}_{n}$ denote the canonical basis of $\mathbb{R}^{n}$. We next introduce the definition of polymatroids, which generalizes matroids.

\begin{definition} \label{def polymatroid}
A \emph{polymatroid} $P=P_{f}\subseteq \mathbb{Z}^{n}$ over $[n]$ with rank function $f$ is  given as $$P=\left\{(a_{1},\ldots,a_{n})\in \mathbb{Z}^{n}\bigg|\sum_{i\in I}a_{i}\leq f(I) \ \text{for any}\ I\subseteq [n] \ \text{and} \ \sum_{i\in [n]}a_{i}=f([n])\right\},$$ where   $f:2^{[n]}\rightarrow \mathbb{Z}$ satisfies

\begin{enumerate}
\item[(i)] $f(\emptyset)=0$;

\item[(ii)] $f(I)+f(J)\geq f(I\cup J)+f(I\cap J)$ for any $I,J\subseteq [n]$ (submodularity).
\end{enumerate}
\end{definition}

A vector $\textbf{a}\in \mathbb{Z}^{n}$ is called a \emph{basis} of $P$ if $\textbf{a}\in P$. It is easy to see that the set of bases (viewed as elements of $\{0,1\}^{n}$) of any matroid is a polymatroid.

In 2022, Bernardi, K\'{a}lm\'{a}n and Postnikov \cite{Bernardi}  defined the polymatroid Tutte polynomial $\mathscr{T}_{P}$ for polymatroids $P$.
\begin{definition} \cite{Bernardi} \label{pop}
Let $P$ be a polymatroid  over $[n]$.
For a basis $\textbf{a}\in P$, an index $i\in [n]$ is \emph{internally active} if $\textbf{a}-\textbf{e}_{i}+\textbf{e}_{j}\notin P$  for any $j<i$. Let $\mathrm{Int}(\textbf{a})=\mathrm{Int}_{P}(\textbf{a})$ denote the set of all
internally active indices with respect to $\textbf{a}$.
An index $i\in [n]$ is \emph{externally active} if $\textbf{a}+\textbf{e}_{i}-\textbf{e}_{j}\notin P$ for any $j<i$. Let $\mathrm{Ext}(\textbf{a})=\mathrm{Ext}_{P}(\textbf{a})$ denote the set of all externally active indices with respect to $\textbf{a}$.

The \emph{polymatroid Tutte polynomial} $\mathscr{T}_{P}(x,y)$ is defined as
$$\mathscr{T}_{P}(x,y):=\sum_{\textbf{a}\in P}x^{oi(\textbf{a})}y^{oe(\textbf{a})}(x+y-1)^{ie(\textbf{a})},$$
where $oi(\textbf{a}):=|\mathrm{Int}(\textbf{a})\setminus \mathrm{Ext}(\textbf{a})|$, $oe(\textbf{a}):=|\mathrm{Ext}(\textbf{a})\setminus \mathrm{Int}(\textbf{a})|$, $ie(\textbf{a}):=|\mathrm{Int}(\textbf{a})\cap \mathrm{Ext}(\textbf{a})|$.
\end{definition}

The order $1<2<\ldots<n$ plays an implicit role in Definition \ref{pop}, but Bernardi et al.~\cite{Bernardi} proved  that $\mathscr{T}_{P}$ depends only on $P$ and not on this order. Similar to the Tutte's original proof for the Tutte polynomial of graphs, the proof in \cite{Guan6} is direct and elementary for it.
\begin{theorem} \label{polymatroid-T} \cite{Bernardi,Guan6}
Let $P$ be a polymatroid. Then $\mathscr{T}_{P}$ only depends on $P$.
\end{theorem}

Obviously, the polymatroid Tutte polynomial generalizes the Tutte polynomial defined by the bases expansion (Definition \ref{Bases}), from matroids to polymatroids. Bernardi et al.~\cite{Bernardi} showed that if $M$ is a matroid of rank $d$ over $[n]$, and  $P=P(M)$ is its corresponding polymatroid, then
$$T_{M}(x,y)=\frac{(x+y-xy)^{n}}{x^{n-d}y^{d}}\mathscr{T}_{P}(\frac{x}{x+y-xy},\frac{y}{x+y-xy}).$$
They \cite{Bernardi} also proved that $\mathscr{T}_{P}(x,y)$ is equivalent to one introduced by Cameron and Fink \cite{Cameron}, which generalizes the Tutte polynomial defined by the rank generating function (Definition \ref{Rank}), from matroids to polymatroids. It is natural to ask the following question.

\begin{question}\label{question}\cite{Bernardi}
Does there exist a more general deletion-contraction relation for the polymatroid Tutte polynomial?
\end{question}

\begin{definition}
Let $P$ be a polymatroid on $[n]$ with rank function $f$. For an index $t\in [n]$, the \emph{deletion} $P\setminus t$ and \emph{contraction} $P/t$, which are polymatroids on $[n]\setminus \{t\}$, are given by the rank functions $f_{P\setminus t}(T)=f(T)$ and $f_{P/t}(T)=f(T\cup \{t\})-f(\{t\})$, for any subset $T\subseteq [n]\setminus \{t\}$, respectively.
 \end{definition}

Let $P$ be a polymatroid over $[n]$ with rank function $f$. For  any $t\in [n]$, let $\alpha_{t}=f([n])-f([n]\setminus \{t\})$, $\beta_{t}=f(\{t\})$ and $T_{t}=\{\alpha_{t}, \alpha_{t}+1,\ldots,\beta_{t}\}$. For any $j\in T_{t}$, define
$P^{t}_{j}:=\{(a_{1},\ldots,a_{n})\in P\mid a_{t}=j\}$ and its projection
\begin{equation*}\label{projection}
\widehat{P}^{t}_{j}:=\{(a_{1},\ldots,a_{t-1},a_{t+1},\ldots,a_{n})\in \mathbb{Z}^{n-1}\mid (a_{1},\ldots,a_{n})\in P^{t}_{j}\}.
\end{equation*}
From \cite{Guan4}, we know that the range $T_{t}$ is chosen such that $P^{t}_{j}$ and $\widehat{P}^{t}_{j}$ are nonempty if and only if $j\in T_{t}$. They also proved that $P^{t}_{j}$ and $\widehat{P}^{t}_{j}$ are polymatroids on $[n]$ and on $[n]\setminus\{t\}$, respectively. Moreover, they found
 \begin{eqnarray*}\label{Del-Con}
 \widehat{P}^{t}_{\alpha_{t}}=P\setminus t\ \text{and} \ \widehat{P}^{t}_{\beta_{t}}=P/t.
\end{eqnarray*}
An recursive formula was obtained for the polymatroid Tutte polynomial as follows in \cite{Bernardi,Guan4}.

\begin{theorem} \label{polymatroid-T} \cite{Bernardi,Guan4}
Let $P$ be a polymatroid over $[n]$. Then for some $t\in [n]$,
\[\mathscr{T}_{P}(x,y)=\left\{\begin{array}{ll}
(x+y-1)\mathscr{T}_{P\setminus t}(x,y),&\text{if}  \ |T_{t}|=1;\\
x\mathscr{T}_{P\setminus t}(x,y)+y\mathscr{T}_{P/t}(x,y)+\sum\limits_{j\in T_{t}\setminus \{\alpha_{t}, \beta_{t}\}}\mathscr{T}_{\widehat{P}^{t}_{j}}(x,y),&\text{if}  \ |T_{t}|\geq 2.
\end{array}\right.\]
\end{theorem}

Clearly, Theorem \ref{polymatroid-T} is consistent with the deletion-contraction formula of the Tutte polynomial of matroids. We next define a polynomial by the recursive formula.

\begin{definition} \label{def-polymatroid-T}
Let $P$ be a polymatroid over $[n]$. A polynomial $\mathscr{T}'_{P}(x,y)$ is defined as follows. If $n=0$, then $\mathscr{T}'_{P}(x,y)=1$. If $n\geq 1$, then for any $t\in [n]$,
\[\mathscr{T}'_{P}(x,y)=\left\{\begin{array}{ll}
(x+y-1)\mathscr{T}'_{P\setminus t}(x,y),&\text{if}  \ |T_{t}|=1;\\
x\mathscr{T}'_{P\setminus t}(x,y)+y\mathscr{T}'_{P/t}(x,y)+\sum\limits_{j\in T_{t}\setminus \{\alpha_{t}, \beta_{t}\}}\mathscr{T}'_{\widehat{P}^{t}_{j}}(x,y),&\text{if}  \ |T_{t}|\geq 2.
\end{array}\right.\]
\end{definition}

In this paper, we shall prove:
\begin{theorem} \label{well-defined}
$\mathscr{T}'_{P}(x,y)$ is well-defined for any polymatroid $P$. Moreover and clearly, $\mathscr{T}'_{P}(x,y)=\mathscr{T}_{P}(x,y)$.
\end{theorem}

\section{Proof of Theorem \ref{well-defined}}
\noindent

We start with two preliminary results.

\begin{proposition}\label{rank function of Pj}\cite{Guan4}
Let $P$ be a polymatroid over $[n]$ with rank function $f$. For some $t\in [n]$, let $f^{t}_{j}$ be the rank function of the polymatroid $\widehat{P}^{t}_{j}$. Then $f^{t}_{j}(I)=\min\{f(I), f(I\cup \{t\})-j\}$ for any subset $I\subseteq [n]\setminus \{t\}$.
\end{proposition}

\begin{proposition}\label{delete} \cite{Guan4}
Let $P$ be a polymatroid  over $[n]$. Then $P\setminus s\setminus t=P\setminus t\setminus s$, $P\setminus s/t=P/t\setminus s$, $P/s\setminus t=P\setminus t/s$ and $P/s/t=P/t/s$ for any $s,t\in [n]$.
\end{proposition}

We now prove Theorem \ref{well-defined}.

\textbf{Proof of Theorem \ref{well-defined}}. Let $P$ be a polymatroid over $[n]$ with rank function $f$. We firstly prove that $\mathscr{T}'_{P}(x,y)$ is processing order independence of the elements in $[n]$ by induction on $n$. If $n=0$,  then $\mathscr{T}'_{P}(x,y)=1$. If $n=1$, then $\mathscr{T}'_{P}(x,y)=x+y-1$. In both cases, the conclusion holds. We now assume that $n\geq 2$.
For any $s,t\in [n]$, we devide into the following three cases to prove it.

\textbf{Case 1.} Suppose that $|T_{t}|=|T_{s}|=1$, that is, $f(\{s\})=f([n])-f([n]\setminus \{s\})$ and $f(\{t\})=f([n])-f([n]\setminus \{t\})$.

Let $f^{t}$ and $f^{s}$ denote the rank functions of polymatroids  $P\setminus t$ and $P\setminus s$. It is easy to see that $f^{t}(\{s\})=f(\{s\})=f([n])-f([n]\setminus \{s\})=(f^{t}([n]\setminus \{t\})+f(\{t\}))-(f^{t}([n]\setminus \{s,t\})+f(\{t\}))=f^{t}([n]\setminus \{t\})-f^{t}([n]\setminus \{s,t\})$ and $f^{s}(\{t\})=f^{s}([n]\setminus \{s\})-f^{s}([n]\setminus \{s,t\})$.

If first deal with $t$, and then deal with $s$, then
\begin{eqnarray*}
	\mathscr{T}'_{P}(x,y)&= & (x+y-1)\mathscr{T}'_{P\setminus t}(x,y)\\
	&=&(x+y-1)^{2}\mathscr{T}'_{P\setminus t\setminus s}(x,y).
\end{eqnarray*}

If first deal with $s$, and then deal with $t$, then
\begin{eqnarray*}
	\mathscr{T}'_{P}(x,y)&=& (x+y-1)\mathscr{T}'_{P\setminus s}(x,y)\\
    &=&(x+y-1)^{2}\mathscr{T}'_{P\setminus s\setminus t}(x,y).
\end{eqnarray*}

By Proposition \ref{delete} and the induction hypothesis, the conclusion is true.

\textbf{Case 2.} Suppose that exactly one of $|T_{t}|$ and $|T_{s}|$ is $1$.

Without loss of generality, we may assume $|T_{t}|=1$ and $|T_{s}|\geq 2$. Let $f^{t}$ and $f^{s}_{i}$ denote the rank functions of polymatroids  $P\setminus t$ and $\widehat{P}^{s}_{i}$ for any $i\in T_{s}$, respectively. It is easy to see that $f^{t}(\{s\})=f(\{s\})=\beta_{s}$, $f^{t}([n]\setminus \{t\})-f^{t}([n]\setminus \{s,t\})=f([n])-f([n]\setminus \{s\})=\alpha_{s}$ and $f^{s}_{i}(\{t\})=f(\{t\})=f([n])-f([n]\setminus \{t\})=(f^{t}([n]\setminus \{s\})+i)-(f^{t}([n]\setminus \{s,t\})+i)=f^{s}_{i}([n]\setminus \{s\})-f^{s}_{i}([n]\setminus \{s,t\})$. Hence, $T^{t}_{s}=\{f^{t}([n]\setminus \{t\})-f^{t}([n]\setminus \{s,t\}),\ldots,f^{t}(\{s\})\}=T_{s}$. Moreover, the following claim holds.

\textbf{Claim 1.} $(\widehat{P}^{s}_{i})\setminus t=(\widehat{P\setminus t})^{s}_{i}$ for any $i\in T_{s}$.

\emph{Proof of Claim 1.} Let   $f_{1}$ and $f_{2}$ be the rank functions of  $(\widehat{P}^{s}_{i})\setminus t$ and $(\widehat{P\setminus t})^{s}_{i}$, respectively.
 Then for any subset $I\subseteq [n]\setminus \{s,t\}$, by Proposition \ref{rank function of Pj},
 $$f_{1}(I)=f^{s}_{i}(I)=\min\{f(I), f(I\cup \{s\})-i\};$$ $$f_{2}(I)=\min\{f^{t}(I), f^{t}(I\cup \{s\})-i\}=\min\{f(I), f(I\cup \{s\})-i\}=f_{1}(I).$$
 Hence, by the definition of the polymatroid, $(\widehat{P}^{s}_{i})\setminus t=(\widehat{P\setminus t})^{s}_{i}$.

If first deal with $t$, and then deal with $s$, then
\begin{eqnarray*}
&& \mathscr{T}'_{P}(x,y)\\
&=& (x+y-1)\mathscr{T}'_{P\setminus t}(x,y)\\
&=&(x+y-1)\left(x\mathscr{T}'_{P\setminus t \setminus s}(x,y)+y\mathscr{T}'_{P\setminus t /s}(x,y)+\sum_{i\in T_{s}\setminus \{\alpha_{s}, \beta_{s}\}}\mathscr{T}'_{(\widehat{P\setminus t})^{s}_{i}}(x,y)\right).
\end{eqnarray*}

If first deal with $s$, and then deal with $t$, then
\begin{eqnarray*}
&& \mathscr{T}'_{P}(x,y)\\
&=& x\mathscr{T}'_{P\setminus s}(x,y)+y\mathscr{T}'_{P/s}(x,y)+\sum_{i\in T_{s}\setminus \{\alpha_{s}, \beta_{s}\}}\mathscr{T}'_{\widehat{P}^{s}_{i}}(x,y)\\
	&=&(x+y-1)x\mathscr{T}'_{P\setminus s\setminus t}(x,y)+(x+y-1)y\mathscr{T}'_{P/s\setminus t}(x,y)\\
	&&+(x+y-1)\left(\sum_{i\in T_{s}\setminus \{\alpha_{s}, \beta_{s}\}}\mathscr{T}'_{(\widehat{P}^{s}_{i})\setminus t}(x,y)\right)\\
	&=&(x+y-1)\left(x\mathscr{T}'_{P\setminus s\setminus t}(x,y)+y\mathscr{T}'_{P/s\setminus t}(x,y)+\sum_{i\in T_{s}\setminus \{\alpha_{s}, \beta_{s}\}}\mathscr{T}'_{(\widehat{P}^{s}_{i})\setminus t}(x,y)\right).
\end{eqnarray*}

By Propositions \ref{delete}, Claim 1, and the induction hypothesis, the conclusion is verified.

\textbf{Case 3.} Suppose that $|T_{t}|\geq 2$ and $|T_{s}|\geq 2$.

Note that  $(P_{i}^{s})_{j}^{t}=(P_{j'}^{t})_{i'}^{s}$ if and only if $i=i'$ and $j=j'$. Moreover, if $(P_{i}^{s})_{j}^{t}=(P_{j'}^{t})_{i'}^{s}$, then $(\widehat{\widehat{P}_{i}^{s}})_{j}^{t}=(\widehat{\widehat{P}_{j'}^{t}})_{i'}^{s}$. But not vice versa. For the sake of brevity, we replace $\textbf{a}\in(P_{i}^{s})_{j}^{t}$ with $\textbf{a}\in (\widehat{\widehat{P}_{i}^{s}})_{j}^{t}$. For any $j\in T_{t}$ and $i\in T_{s}$, let $f^{t}_{j}$ and $f^{s}_{i}$ denote the rank functions of polymatroids  $\widehat{P}^{t}_{i}$ and $\widehat{P}^{s}_{j}$, respectively. Denote $\alpha^{tj}_{s}=f^{t}_{j}([n]\setminus\{t\})-f^{t}_{j}([n]\setminus\{s,t\})$, $\beta^{tj}_{s}=f^{t}_{j}(\{s\})$, $T^{tj}_{s}=\{\alpha^{tj}_{s},\alpha^{tj}_{s}+1,\ldots,\beta^{tj}_{s}\}$, $\alpha^{si}_{t}=f^{s}_{i}([n]\setminus\{s\})-f^{s}_{i}([n]\setminus\{s,t\})$, $\beta^{si}_{t}=f^{s}_{i}(\{t\})$ and $T^{si}_{t}=\{\alpha^{si}_{t},\alpha^{si}_{t}+1,\ldots,\beta^{si}_{t}\}$. Obviously, $|T^{tj}_{s}|\geq 1$ for any $j\in T_{t}$, and $|T^{si}_{t}|\geq 1$ for any $i\in T_{s}$. By Proposition \ref{rank function of Pj}, we have that
$$\alpha^{t\alpha_{t}}_{s}=f([n]\setminus\{t\})-f([n]\setminus\{s,t\});$$
$$\beta^{t\alpha_{t}}_{s}=f(\{s\})=\beta_{s};$$
 \begin{eqnarray*}
	\alpha^{t\beta_{t}}_{s}&= & (f([n])-f(\{t\}))-(f([n]\setminus\{s\})-f(\{t\}))\\
	&=&f([n])-f([n]\setminus\{s\})=\alpha_{s};
\end{eqnarray*}
$$\beta^{t\beta_{t}}_{s}=f(\{s,t\})-f(\{t\}).$$
Similarly, we have
$$\alpha^{s\alpha_{s}}_{t}=f([n]\setminus\{s\})-f([n]\setminus\{s,t\};$$
$$\beta^{s\alpha_{s}}_{t}=f(\{t\})=\beta_{t};$$
$$\alpha^{s\beta_{s}}_{t}=f([n])-f([n]\setminus\{t\}=\alpha_{t};$$
$$\beta^{s\beta_{s}}_{t}=f(\{s,t\})-f(\{s\}).$$

 We next devide into two sub-cases to prove it.

\textbf{Sub-Case 3.1.} Assume that $|T^{tj}_{s}|\geq 2$ for any $j\in T_{t}$, and $|T^{si}_{t}|\geq 2$ for any $i\in T_{s}$.

If first deal with $t$, and then deal with $s$, then
\begin{eqnarray*}
&& \mathscr{T}'_{P}(x,y)\\
&=& x\mathscr{T}'_{P\setminus t}(x,y)+y\mathscr{T}'_{P/t}(x,y)+\sum_{j\in T_{t}\setminus \{\alpha_{t}, \beta_{t}\}}\mathscr{T}'_{\widehat{P}^{t}_{j}}(x,y)\\
	&= &x\left(x\mathscr{T}'_{P\setminus t \setminus s}(x,y)+y\mathscr{T}'_{P\setminus t /s}(x,y)+\sum_{i\in T^{t\alpha_{t}}_{s}\setminus \{\alpha^{t\alpha_{t}}_{s}, \beta^{t\alpha_{t}}_{s}\}}\mathscr{T}'_{(\widehat{P\setminus t})^{s}_{i}}(x,y)\right)\\
& &+y\left(x\mathscr{T}'_{P/t\setminus s}(x,y)+y\mathscr{T}'_{P/t/s}(x,y)+\sum_{i\in T^{t\beta_{t}}_{s}\setminus \{\alpha^{t\beta_{t}}_{s}, \beta^{t\beta_{t}}_{s}\}}\mathscr{T}'_{(\widehat{P/t})^{s}_{i}}(x,y)\right)\\
& &+\sum_{j\in T_{t}\setminus \{\alpha_{t}, \beta_{t}\}}\left(x\mathscr{T}'_{(\widehat{P}^{t}_{j})\setminus s}(x,y)+y\mathscr{T}'_{(\widehat{P}^{t}_{j}) /s}(x,y)+\sum_{i\in T^{tj}_{s}\setminus \{\alpha^{tj}_{s}, \beta^{tj}_{s}\}}\mathscr{T}'_{(\widehat{\widehat{P}^{t}_{j}})^{s}_{i}}(x,y)\right)\\
&= &x^{2}\mathscr{T}'_{P\setminus t \setminus s}(x,y)+y^{2}\mathscr{T}'_{P/t/s}(x,y)+xy(\mathscr{T}'_{P\setminus t /s}(x,y)+\mathscr{T}'_{P/t\setminus s}(x,y))\\
& &+x\left(\sum_{i\in T^{t\alpha_{t}}_{s}\setminus \{\alpha^{t\alpha_{t}}_{s}, \beta^{t\alpha_{t}}_{s}\}}\mathscr{T}'_{(\widehat{P\setminus t})^{s}_{i}}(x,y)+\sum_{j\in T_{t}\setminus \{\alpha_{t}, \beta_{t}\}}\mathscr{T}'_{(\widehat{P}^{t}_{j})\setminus s}(x,y)\right)\\
& &+y\left(\sum_{i\in T^{t\beta_{t}}_{s}\setminus \{\alpha^{t\beta_{t}}_{s}, \beta^{t\beta_{t}}_{s}\}}\mathscr{T}'_{(\widehat{P/t})^{s}_{i}}(x,y)+\sum_{j\in T_{t}\setminus \{\alpha_{t}, \beta_{t}\}}\mathscr{T}'_{(\widehat{P}^{t}_{j})/s}(x,y)\right)\\
& &+\sum_{j\in T_{t}\setminus \{\alpha_{t}, \beta_{t}\}}\sum_{i\in T^{tj}_{s}\setminus \{\alpha^{tj}_{s}, \beta^{tj}_{s}\}}\mathscr{T}'_{(\widehat{\widehat{P}^{t}_{j}})^{s}_{i}}(x,y).
\end{eqnarray*}

If first deal with $s$, and then deal with $t$, then
\begin{eqnarray*}
&& \mathscr{T}'_{P}(x,y)\\
&=& x^{2}\mathscr{T}'_{P\setminus s \setminus t}(x,y)+y^{2}\mathscr{T}'_{P/s/t}(x,y)+xy(\mathscr{T}'_{P\setminus s /t}(x,y)+\mathscr{T}'_{P/s\setminus t}(x,y))\\
& &+x\left(\sum_{j\in T^{s\alpha_{s}}_{t}\setminus \{\alpha^{s\alpha_{s}}_{t}, \beta^{s\alpha_{s}}_{t}\}}\mathscr{T}'_{(\widehat{P\setminus s})^{t}_{j}}(x,y)+\sum_{i\in T_{s}\setminus \{\alpha_{s}, \beta_{s}\}}\mathscr{T}'_{(\widehat{P}^{s}_{i})\setminus t}(x,y)\right)\\
& &+y\left(\sum_{j\in T^{s\beta_{s}}_{t}\setminus \{\alpha^{s\beta_{s}}_{t}, \beta^{s\beta_{s}}_{t}\}}\mathscr{T}'_{(\widehat{P/s})^{t}_{j}}(x,y)+\sum_{i\in T_{s}\setminus \{\alpha_{s}, \beta_{s}\}}\mathscr{T}'_{(\widehat{P}^{s}_{i})/t}(x,y)\right)\\
& &+\sum_{i\in T_{s}\setminus \{\alpha_{s}, \beta_{s}\}}\sum_{j\in T^{si}_{t}\setminus \{\alpha^{si}_{t}, \beta^{si}_{t}\}}\mathscr{T}'_{(\widehat{\widehat{P}^{s}_{i}})^{t}_{j}}(x,y).
\end{eqnarray*}
By Proposition \ref{delete} and the induction hypothesis,  coefficients of $x^{2}$, $y^{2}$ and $xy$ are same.
Note that,
\begin{enumerate}
\item[(I-i)] $a_{t}=\alpha_{t}$ and $a_{s}=\alpha^{t\alpha_{t}}_{s}=f([n]\setminus \{t\})-f([n]\setminus \{s,t\})$ for any $\textbf{a}\in P\setminus t\setminus s$;

\item[(I-ii)]  $a_{s}=\alpha_{s}$ and $a_{t}=\alpha^{s\alpha_{s}}_{t}=f([n]\setminus \{s\})-f([n]\setminus \{s,t\})$ for any $\textbf{a}\in P\setminus s\setminus t$;

\item[(I-iii)] $a_{t}=\alpha_{t}$ and $a_{s}=i$ for any $i\in T^{t\alpha_{t}}_{s}\setminus \{\alpha^{t\alpha_{t}}_{s}, \beta^{t\alpha_{t}}_{s}\}$ and for any $\textbf{a}\in (\widehat{P\setminus t})^{s}_{i}$;

\item[(I-iv)] for any $j\in T_{t}\setminus \{\alpha_{t}, \beta_{t}\}$ and for any $\textbf{a}\in(\widehat{P}^{t}_{j})\setminus s$, we have $a_{t}=j$. Since $j>\alpha_{t}=f([n])-f([n]\setminus \{t\})$, we have that
\begin{eqnarray*}
	a_{s}&= & f^{t}_{j}([n]\setminus \{t\})-f^{t}_{j}([n]\setminus \{s,t\})\\
&=&\min\{f([n]\setminus \{t\}), f([n])-j\}-\min\{f([n]\setminus \{s,t\}), f([n]\setminus \{s\})-j\}\\
&=&f([n])-j-\min\{f([n]\setminus \{s,t\}), f([n]\setminus \{s\})-j\}.
\end{eqnarray*}
More precisely,
\[a_{s}=\left\{\begin{array}{ll}
f([n])-f([n]\setminus \{s\})=\alpha_{s},&\text{if}  \ \alpha^{s\alpha_{s}}_{t}\leq j< f(\{t\})=\beta^{s\alpha_{s}}_{t}=\beta_{t},\\
f([n])-j-f([n]\setminus \{s,t\}),&\text{if}  \ \alpha_{t}< j<\alpha^{s\alpha_{s}}_{t};
\end{array}\right.\]
(That is, if $\alpha_{s}=f([n])-f([n]\setminus \{s\})< a_{s}< f([n]\setminus \{t\})-f([n]\setminus \{s,t\})=\alpha^{t\alpha_{t}}_{s}$, then $a_{s}+a_{t}=f([n])-f([n]\setminus \{s,t\})$.)
\item[(I-v)] for any $j\in T^{s\alpha_{s}}_{t}\setminus \{\alpha^{s\alpha_{s}}_{t}, \beta^{s\alpha_{s}}_{t}\}$ and for any $\textbf{a}\in (\widehat{P\setminus s})^{t}_{j}$, we have $a_{s}=\alpha_{s}$ and $a_{t}=j$;

\item[(I-vi)] for any $i\in T_{s}\setminus \{\alpha_{s}, \beta_{s}\}$ and for any $\textbf{a}\in(\widehat{P}^{s}_{i})\setminus t$, we have $a_{s}=i$,  and
\[a_{t}=\left\{\begin{array}{ll}
f([n])-f([n]\setminus \{t\})=\alpha_{t},&\text{if}  \ \alpha^{t\alpha_{t}}_{s}\leq i< f(\{s\})=\beta^{t\alpha_{t}}_{s}=\beta_{s},\\
f([n])-i-f([n]\setminus \{s,t\}),&\text{if}  \ \alpha_{s}< i< \alpha^{t\alpha_{t}}_{s}.
\end{array}\right.\]
(In this case, $a_{s}+a_{t}=f([n])-f([n]\setminus \{s,t\})$ and $\alpha_{s}< a_{s}< \alpha^{t\alpha_{t}}_{s}$.)
\end{enumerate}

Hence, $(\widehat{P\setminus t})^{s}_{i}=(\widehat{P}^{s}_{i})\setminus t$ for any $i\in T^{t\alpha_{t}}_{s}\setminus \{\alpha^{t\alpha_{t}}_{s}, \beta^{t\alpha_{t}}_{s}\}$,  $(\widehat{P\setminus s})^{t}_{j}=(\widehat{P}^{t}_{j})\setminus s$ for any $j\in T^{s\alpha_{s}}_{t}\setminus \{\alpha^{s\alpha_{s}}_{t}\beta^{s\alpha_{s}}_{t}\}$, $(\widehat{P}^{t}_{j})\setminus s=(\widehat{P}^{s}_{i})\setminus t$ for any $\alpha_{t}< j<\alpha^{s\alpha_{s}}_{t}$ and $i=f([n])-j-f([n]\setminus \{s,t\})$, $(\widehat{P}^{t}_{j})\setminus s=P\setminus s\setminus t=P\setminus t\setminus s=(\widehat{P}^{s}_{i})\setminus t$ when $j= f([n]\setminus\{s\})-f([n]\setminus \{s,t\})$ and $i=f([n]\setminus\{t\})-f([n]\setminus \{s,t\})$. By induction hypothesis,
\begin{eqnarray*}
&& \sum_{i\in T^{t\alpha_{t}}_{s}\setminus \{\alpha^{t\alpha_{t}}_{s}, \beta^{t\alpha_{t}}_{s}\}}\mathscr{T}'_{(\widehat{P\setminus t})^{s}_{i}}(x,y)+\sum_{j\in T_{t}\setminus \{\alpha_{t}, \beta_{t}\}}\mathscr{T}'_{(\widehat{P}^{t}_{j})\setminus s}(x,y)\\
&=& \sum_{i\in T^{t\alpha_{t}}_{s}\setminus \{\alpha^{t\alpha_{t}}_{s}, \beta^{t\alpha_{t}}_{s}\}}\mathscr{T}'_{(\widehat{P\setminus t})^{s}_{i}}(x,y)+\sum_{\alpha_{t}< j<\alpha^{s\alpha_{s}}_{t}}\mathscr{T}'_{(\widehat{P}^{t}_{j})\setminus s}(x,y)\\
&+&\mathscr{T}'_{(\widehat{P}^{t}_{\alpha^{s\alpha_{s}}_{t}})\setminus s}(x,y)+\sum_{j\in T^{s\alpha_{s}}_{t}\setminus \{\alpha^{s\alpha_{s}}_{t}\beta^{s\alpha_{s}}_{t}\}}\mathscr{T}'_{(\widehat{P}^{t}_{j})\setminus s}(x,y)\\
	&=&\sum_{i\in T^{t\alpha_{t}}_{s}\setminus \{\alpha^{t\alpha_{t}}_{s}, \beta^{t\alpha_{t}}_{s}\}}\mathscr{T}'_{(\widehat{P}^{s}_{i})\setminus t}(x,y)+\sum_{\alpha_{s}< i< \alpha^{t\alpha_{t}}_{s}}\mathscr{T}'_{(\widehat{P}^{s}_{i})\setminus t}(x,y)\\
	&+&\mathscr{T}'_{(\widehat{P}^{s}_{\alpha^{t\alpha_{t}}_{s}})\setminus t}(x,y)+\sum_{j\in T^{s\alpha_{s}}_{t}\setminus \{\alpha^{s\alpha_{s}}_{t}, \beta^{s\alpha_{s}}_{t}\}}\mathscr{T}'_{(\widehat{P\setminus s})^{t}_{j}}(x,y)\\
&=&\sum_{j\in T^{s\alpha_{s}}_{t}\setminus \{\alpha^{s\alpha_{s}}_{t}, \beta^{s\alpha_{s}}_{t}\}}\mathscr{T}'_{(\widehat{P\setminus s})^{t}_{j}}(x,y)+\sum_{i\in T_{s}\setminus \{\alpha_{s}, \beta_{s}\}}\mathscr{T}'_{(\widehat{P}^{s}_{i})\setminus t}(x,y).
\end{eqnarray*}
Hence, the coefficients of $x$ are same.

Similarly,
\begin{enumerate}
\item[(II-i)] $a_{t}=\beta_{t}$ and $a_{s}=f(\{s,t\})-\beta_{t}=f(\{s,t\})-f(\{t\})$ for any $\textbf{a}\in P/t/s$;

\item[(II-ii)]  $a_{s}=\beta_{s}$ and $a_{t}=f(\{s,t\})-\beta_{s}=f(\{s,t\})-f(\{s\})$ for any $\textbf{a}\in P/s/t$;

\item[(II-iii)] $a_{t}=\beta_{t}$ and $a_{s}=i$ for any $i\in T^{t\beta_{t}}_{s}\setminus \{\alpha^{t\beta_{t}}_{s}, \beta^{t\beta_{t}}_{s}\}$ and for any $\textbf{a}\in(\widehat{P/t})^{s}_{i}$;

\item[(II-iv)] for any $j\in T_{t}\setminus \{\alpha_{t}, \beta_{t}\}$, and for any $\textbf{a}\in(\widehat{P}^{t}_{j})/s$, we have $a_{t}=j$ and ~$a_{s}=f^{t}_{j}(\{s\})=\min\{f(\{s\}), f(\{s,t\})-j\}$.
More precisely,
\[a_{s}=\left\{\begin{array}{ll}
f(\{s\})=\beta_{s},&\text{if}  \ \alpha^{s\beta_{s}}_{t}=\alpha_{t}< i\leq \beta^{s\beta_{s}}_{t},\\
f(\{s,t\})-j,&\text{if}  \ \beta^{s\beta_{s}}_{t}<j < f(\{t\})=\beta_{t};
\end{array}\right.\]
(In this case, $a_{s}+a_{t}=f(\{s,t\})$ and $\beta^{t\beta_{t}}_{s}=f(\{s,t\})-f(\{t\})< a_{s}< f(\{s\})=\beta_{s}$ if $\beta^{s\beta_{s}}_{t}<j < \beta_{t}$.)
\item[(II-v)] for any $j\in T^{s\beta_{s}}_{t}\setminus \{\alpha^{s\beta_{s}}_{t}, \beta^{s\beta_{s}}_{t}\}$ and for any $\textbf{a}\in(\widehat{P/s})^{t}_{j}$, we have  $a_{s}=\beta_{s}$ and $a_{t}=j$;
\item[(II-vi)] for any $i\in T_{s}\setminus \{\alpha_{s}, \beta_{s}\}$ and for any $\textbf{a}\in(\widehat{P}^{s}_{i})/t$, we have $a_{s}=i$ and ~$a_{t}=f^{s}_{i}(\{t\})=\min\{f(\{t\}), f(\{s,t\})-i\}$.
More precisely,
\[a_{t}=\left\{\begin{array}{ll}
f(\{t\})=\beta_{t},&\text{if}  \ \alpha^{t\beta_{t}}_{s}=\alpha_{s}< j\leq \beta^{t\beta_{t}}_{s},\\
f(\{s,t\})-i,&\text{if}  \ \beta^{t\beta_{t}}_{s}<i < f(\{s\})=\beta_{s}.
\end{array}\right.\]
(In this case, $a_{s}+a_{t}=f(\{s,t\})$ and $\beta^{s\beta_{s}}_{t}<a_{t} < \beta_{t}$ when $\beta^{t\beta_{t}}_{s}< i< \beta_{s}$.)
\end{enumerate}

Hence, $(\widehat{P/t})^{s}_{i}=(\widehat{P}^{s}_{i})/t$ for any $i\in T^{t\beta_{t}}_{s}\setminus \{\alpha^{t\beta_{t}}_{s}, \beta^{t\beta_{t}}_{s}\}$, $(\widehat{P/s})^{t}_{j}=(\widehat{P}^{t}_{j})/s$ for any $j\in T^{s\beta_{s}}_{t}\setminus \{\alpha^{s\beta_{s}}_{t}, \beta^{s\beta_{s}}_{t}\}$, $(\widehat{P}^{t}_{j})/s=(\widehat{P}^{s}_{i})/t$ for any $f(\{s,t\})-f(\{s\})<j < f(\{t\})$ and $i=f(\{s,t\})-j$, and $(\widehat{P}^{t}_{j})/s=P/s/t=P/t/s=(\widehat{P}^{s}_{i})/t$ when $j= f(\{s,t\})-f(\{s\})$ and $i=f(\{s,t\})$. By induction hypothesis,  coefficients of $y$ are same.

 Moreover, we have that
 \begin{enumerate}
 \item[(III-i)] $a_{t}=\alpha_{t}$ and~$a_{s}=\beta^{t\alpha_{t}}_{s}=\beta_{s}$ for any $\textbf{a}\in P\setminus t /s$;
 \item[(III-ii)] $a_{t}=\beta_{t}$ and~$a_{s}=\alpha_{s}=\alpha^{t\beta_{t}}_{s}$ for any $\textbf{a}\in P/t\setminus s$;
 \item[(III-iii)] $a_{s}=\alpha_{s}$ and~$a_{t}=\beta^{s\alpha_{s}}_{t}=\beta_{t}$ for any~$\textbf{a}\in P\setminus s /t$;
 \item[(III-iv)] $a_{s}=\beta_{s}$ and~$a_{t}=\alpha_{t}=\alpha^{s\beta_{s}}_{t}$ for any~$\textbf{a}\in P/s\setminus t$.
\end{enumerate}

According to the above discussion (see Table 1), corresponding polymatroids of the constant term are same. Therefore, by induction hypothesis, their corresponding coefficients of the constant term are same.

\begin{table}[htp]
\begin{center}
\renewcommand\arraystretch{1.2}\footnotesize
\caption{A comparison of the polymatroid Tutte polynomial obtained by dealing with $s,t$.}
 \begin{tabular}{|c |c|c|c|c|c|c|c|c|c|c |c|c|c|c|}
  \hline
 \text{first}~$t$  \text{and then}~$s$ &$a_{t}$& $a_{s}$& \text{first}~$s$  \text{and then}~$t$\\
 \hline
 $\textbf{a}\in P\setminus t\setminus s$& $\alpha_{t}$ & $f([n]\setminus \{t\})-$ & $\textbf{a}\in(\widehat{P}^{s}_{i})\setminus t$, \text{where}~$i=$\\
 & &$f([n]\setminus \{s,t\})$ &$f([n]\setminus\{t\})-f([n]\setminus \{s,t\})$\\
 \hline
$\textbf{a}\in P\setminus t /s$ & $\alpha_{t}$ & $\beta_{s}$ & $\textbf{a}\in P/s\setminus t$\\
 \hline
 $\textbf{a}\in P/t/s$&$\beta_{t}$& $f(\{s,t\})-$& $\textbf{a}\in(\widehat{P}^{s}_{i})/ t$, \text{where}\\
 & &$f(\{t\})$ &$i=f(\{s,t\})-f(\{t\})$\\
 \hline
 $\textbf{a}\in P/t\setminus s$&$\beta_{t}$& $\alpha_{s}$& $\textbf{a}\in P\setminus s /t$\\
 \hline
 $\textbf{a}\in (\widehat{P\setminus t})^{s}_{i}$, \text{where}, &$\alpha_{t}$& $i$& $\textbf{a}\in(\widehat{P}^{s}_{i})\setminus t$, \text{where}~ \\

$i\in T^{t\alpha_{t}}_{s}\setminus \{\alpha^{t\alpha_{t}}_{s}, \beta^{t\alpha_{t}}_{s}\}$& & &$i\in T^{t\alpha_{t}}_{s}\setminus \{\alpha^{t\alpha_{t}}_{s}, \beta^{t\alpha_{t}}_{s}\}$\\
 \hline
 $\textbf{a}\in(\widehat{P}^{t}_{j})\setminus s$, \text{where}&$j$& $\alpha_{s}$& $\textbf{a}\in (\widehat{P\setminus s})^{t}_{j}$, \text{where}\\

$j\in T^{s\alpha_{s}}_{t}\setminus \{\alpha^{s\alpha_{s}}_{t}, \beta^{s\alpha_{s}}_{t}\}$&& &$j\in T^{s\alpha_{s}}_{t}\setminus \{\alpha^{s\alpha_{s}}_{t}, \beta^{s\alpha_{s}}_{t}\}$ \\
 \hline
 $\textbf{a}\in(\widehat{P}^{t}_{j})\setminus s$,~\text{where}~$j=$&$f([n]\setminus \{s\})$& $\alpha_{s}$& $\textbf{a}\in P\setminus s\setminus t$\\

$f([n]\setminus\{s\})-f([n]\setminus \{s,t\})$&$-f([n]\setminus \{s,t\})$& & \\
 \hline
 $\textbf{a}\in(\widehat{P}^{t}_{j})\setminus s$, \text{where}~$\alpha_{t}<j<$&$j$& $f([n])-j-$& $\textbf{a}\in(\widehat{P}^{s}_{i})\setminus t$, \text{where}~$i=$\\
$f([n]\setminus\{s\})-f([n]\setminus \{s,t\})$& &$f([n]\setminus \{s,t\})$ & $f([n])-j-f([n]\setminus \{s,t\})$\\
 \hline
 $\textbf{a}\in (\widehat{P/t})^{s}_{i}$, \text{where}&$\beta_{t}$& $i$& $\textbf{a}\in(\widehat{P}^{s}_{i})/t$, \text{where}~ \\

$i\in T^{t\beta_{t}}_{s}\setminus \{\alpha^{t\beta_{t}}_{s}, \beta^{t\beta_{t}}_{s}\}$& & &$i\in T^{t\beta_{t}}_{s}\setminus \{\alpha^{t\beta_{t}}_{s}, \beta^{t\beta_{t}}_{s}\}$\\
 \hline
 $\textbf{a}\in(\widehat{P}^{t}_{j})/s$, \text{where}&$j$& $\beta_{s}$& $\textbf{a}\in (\widehat{P/s})^{t}_{j}$, \text{where}\\

$j\in T^{s\beta_{s}}_{t}\setminus \{\alpha^{s\beta_{s}}_{t}, \beta^{s\beta_{s}}_{t}\}$&& &$j\in T^{s\beta_{s}}_{t}\setminus \{\alpha^{s\beta_{s}}_{t}, \beta^{s\beta_{s}}_{t}\}$ \\
 \hline
 $\textbf{a}\in(\widehat{P}^{t}_{j})/s$,~\text{where}~&$f(\{s,t\})$& $\beta_{s}$& $\textbf{a}\in P/s/t$\\

$j=f(\{s,t\})-f(\{s\})$&$-f(\{s\})$& & \\
 \hline
 $\textbf{a}\in(\widehat{P}^{t}_{j})/s$, \text{where}&$j$& $f(\{s,t\})-j$& $\textbf{a}\in(\widehat{P}^{s}_{i})/t$, \text{where}\\
$f(\{s,t\})-f(\{s\})<j<\beta_{t}$& && $i=f(\{s,t\})-j$\\
 \hline
\end{tabular}
\end{center}
\end{table}

\textbf{Sub-Case 3.2.}  $|T^{tj}_{s}|=1$ for some $j\in T_{t}$ or  $|T^{si}_{t}|=1$ for some $i\in T_{s}$.

In this subcase, the following claims hold.

\textbf{Claim 2.}  $|T^{tj}_{s}|=1$ for some $j\in T_{t}$ if and only if one of the following holds:

\begin{enumerate}
\item[(i)] $f([n]\setminus \{s\})+f(\{s,t\})= f([n])+f(\{t\})$ and~$j=f(\{t\})=\beta_{t}$;
\item[(ii)] $f(\{s,t\})+f([n]\setminus \{s,t\})=f([n])$ and~$j$ is~$T_{t}$~an arbitrary integer;
\item[(iii)] $f([n]\setminus \{s,t\})+f(\{s\})= f([n]\setminus \{t\})$ and~$j=\alpha_{t}$.
\end{enumerate}

\emph{Sufficiency.} (i) If $j=f(\{t\})=\beta_{t}$, then by the submodularity of $f$, we have that $f^{t}_{j}(\{s\})=f(\{s,t\})-f(\{t\})$,  $f^{t}_{j}([n]\setminus\{t\})=f([n])-f(\{t\})$ and ~$f^{t}_{j}([n]\setminus\{s,t\})=f([n]\setminus\{s\})-f(\{t\})$. Note that $f([n]\setminus \{s\})+f(\{s,t\})= f([n])+f(\{t\})$. So, $f^{t}_{j}(\{s\})=f^{t}_{j}([n]\setminus\{t\})-f^{t}_{j}([n]\setminus\{s,t\})$, that is, $|T^{tj}_{s}|=1$.

(ii) If $f(\{s,t\})+f([n]\setminus \{s,t\})=f([n])$, then $a_{s}+a_{t}=f(\{s,t\})$ for any~$\textbf{a}\in P$, that is, $a_{t}=f(\{s,t\})-a_{s}$. Hence, $f^{t}_{j}(\{s\})=f^{t}_{j}([n]\setminus\{t\})-f^{t}_{j}([n]\setminus\{s,t\})$.

(iii) If $j=f([n])-f([n]\setminus\{t\})=\alpha_{t}$, then by the submodularity of $f$, we have that $f^{t}_{j}(\{s\})=f(\{s\})$,  $f^{t}_{j}([n]\setminus\{t\})=f([n]\setminus\{t\})$ and~$f^{t}_{j}([n]\setminus\{s,t\})=f([n]\setminus\{s,t\})$.  Note that $f([n]\setminus \{s,t\})+f(\{s\})= f([n]\setminus \{t\})$. Hence, $f^{t}_{j}(\{s\})=f^{t}_{j}([n]\setminus\{t\})-f^{t}_{j}([n]\setminus\{s,t\})$.

\emph{Necessity.}  Note that $f_{j}^{t}(s)=\min\{f(\{s\}),f(\{s,t\})-j\}$, $f_{j}^{t}([n]\setminus \{s,t\})=\min\{f([n]\setminus \{s,t\}),f([n]\setminus \{s\})-j\}$, $f_{j}^{t}([n]\setminus \{t\})=f([n])-j$ by Proposition \ref{rank function of Pj} and $f_{j}^{t}(s)+f_{j}^{t}([n]\setminus \{s,t\})=f_{j}^{t}([n]\setminus \{t\})$ by $|T^{tj}_{s}|=1$.

If $j\geq f(\{s,t\})-f(\{s\})$ and $j\geq f([n]\setminus \{s\})-f([n]\setminus \{s,t\})$, then $f(\{s,t\})-j=f([n])-f([n]\setminus \{s\})$, that is, $j=f(\{s,t\})-(f([n])-f([n]\setminus \{s\}))$. By the submodularity of $f$, we have that $f([n]\setminus \{s\})+f(\{s,t\})\geq f([n])+f(\{t\})$. Hence, $j\geq f(\{t\})$. Since $j\leq f(\{t\})$, we have $j=f(\{t\})=\beta_{t}$.

If $f([n]\setminus \{s\})-f([n]\setminus \{s,t\})\leq j\leq f(\{s,t\})-f(\{s\})$, then $f(\{s\})=f([n])-f([n]\setminus \{s\})$, that is, $|T_{s}|=1$, a contradiction.

If~$f(\{s,t\})-f(\{s\})\leq j\leq f([n]\setminus \{s\})-f([n]\setminus \{s,t\})$, then $f(\{s,t\})+f([n]\setminus \{s,t\})=f([n])$. Note that in this case,~$f(\{s,t\})-f(\{s\})=f([n])-f([n]\setminus \{s,t\})-f(\{s\})\leq f([n])-f([n]\setminus \{t\})=\alpha_{t}$ and $f([n]\setminus \{s\})-f([n]\setminus \{s,t\})=f([n]\setminus \{s\})-(f([n])-f(\{s,t\}))\geq f(\{t\})=\beta_{t}$, Hence, $j$ is an arbitrary integer in $T_{t}$.

If~$j\leq f([n]\setminus \{s\})-f([n]\setminus \{s,t\})$ and~$j\leq f(\{s,t\})-f(\{s\})$, then $j=f([n])-f([n]\setminus \{s,t\})-f(\{s\})$. By the submodularity of $f$, we have that $f([n]\setminus \{s,t\})+f(\{s\})\geq f([n]\setminus \{t\})$. Hence, $j\leq f([n])-f([n]\setminus \{t\})=\alpha_{t}$. Since~$j\geq \alpha_{t}$, we have $j=\alpha_{t}$.

Moreover, we have the following result similar to Claim 1.

\textbf{Claim 3.}  $|T^{si}_{t}|=1$ for some $i\in T_{s}$ if and only if one of the following holds:
\begin{enumerate}
\item[(i)] $f([n]\setminus \{t\})+f(\{s,t\})= f([n])+f(\{s\})$ and $i=f(\{s\})=\beta_{s}$;
\item[(ii)] $f(\{s,t\})+f([n]\setminus \{s,t\})=f([n])$ and $i$ is an arbitrary integer in $T_{s}$;
\item[(iii)] $f([n]\setminus \{s,t\})+f(\{t\})= f([n]\setminus \{s\})$ and~$i=\alpha_{s}$.
\end{enumerate}

\textbf{(3-2-1)} Assume that~$f(\{s,t\})+f([n]\setminus \{s,t\})=f([n])$.

By Claims 2 and 3, we have $|T^{tj}_{s}|=1$ for any $j\in T_{t}$, and $|T^{si}_{t}|=1$ for any $i\in T_{s}$. Moreover, $|T_{t}|=|T_{s}|$ and $P\setminus t\setminus s=P\setminus s\setminus t=P/t\setminus s=P/s\setminus t=(\widehat{P}^{t}_{j})\setminus s=(\widehat{P}^{s}_{i})\setminus t$ for any $j\in T_{t}$ and for any $i\in T_{s}$.

If first deal with $t$, and then deal with $s$, then
 \begin{eqnarray*}
 && \mathscr{T}'_{P}(x,y)\\
&=& x\mathscr{T}'_{P\setminus t}(x,y)+y\mathscr{T}'_{P/t}(x,y)+\sum_{j\in T_{t}\setminus \{\alpha_{t}, \beta_{t}\}}\mathscr{T}'_{\widehat{P}^{t}_{j}}(x,y)\\
	&= &(x+y-1)\left(x\mathscr{T}'_{P\setminus t\setminus s}(x,y)+y\mathscr{T}'_{P/t\setminus s}(x,y)+\sum_{j\in T_{t}\setminus \{\alpha_{t}, \beta_{t}\}}\mathscr{T}'_{(\widehat{P}^{t}_{j})\setminus s}(x,y)\right)\\
	&= &(x+y-1)(x+y+|T_{t}|-2)(\mathscr{T}'_{P\setminus t\setminus s}(x,y)).
\end{eqnarray*}

If first deal with $s$, and then deal with $t$, then $$\mathscr{T}'_{P}(x,y)=(x+y-1)(x+y+|T_{s}|-2)\mathscr{T}'_{P\setminus s\setminus t}(x,y).$$

Therefore, by induction hypothesis and Proposition \ref{delete}, the conclusion holds.

We now assume $f(\{s,t\})+f([n]\setminus \{s,t\})>f([n])$. It is easy to see that $f([n]\setminus \{s\})+f(\{s,t\})= f([n])+f(\{t\})$ and $f([n]\setminus \{s,t\})+f(\{t\})=f([n]\setminus \{s\})$ do not hold simultaneously, and $f([n]\setminus \{s,t\})+f(\{s\})= f([n]\setminus \{t\})$ and~$f([n]\setminus \{t\})+f(\{s,t\})=f([n])+f(\{s\})$ do not hold simultaneously. Hence, it is enough to consider the following cases.

\textbf{(3-2-2)} Assume that $f([n]\setminus \{s\})+f(\{s,t\})= f([n])+f(\{t\})$ and $f([n]\setminus \{s,t\})+f(\{s\})= f([n]\setminus \{t\})$.
Then $\alpha^{t\beta_{t}}_{s}=\beta^{t\beta_{t}}_{s}$, $\alpha^{t\alpha_{t}}_{s}=\beta^{t\alpha_{t}}_{s}$, $\alpha^{s\alpha_{s}}_{t}<\beta^{s\alpha_{s}}_{t}$ and $\alpha^{s\beta_{s}}_{t}<\beta^{s\beta_{s}}_{t}$. Moreover, $P\setminus t\setminus s=P\setminus t/s$ and $P/ t\setminus s=P/ t/s$.

If first deal with $t$, and then deal with $s$, then
\begin{eqnarray*}
&&	\mathscr{T}'_{P}(x,y)\\
&=& x\mathscr{T}'_{P\setminus t}(x,y)+y\mathscr{T}'_{P/t}(x,y)+\sum_{j\in T_{t}\setminus \{\alpha_{t}, \beta_{t}\}}\mathscr{T}'_{\widehat{P}^{t}_{j}}(x,y)\\
	&= &x(x+y-1)\mathscr{T}'_{P\setminus t \setminus s}(x,y)+y(x+y-1)\mathscr{T}'_{P/t\setminus s}(x,y)\\
& &+\sum_{j\in T_{t}\setminus \{\alpha_{t}, \beta_{t}\}}\left(x\mathscr{T}'_{(\widehat{P}^{t}_{j})\setminus s}(x,y)+y\mathscr{T}'_{(\widehat{P}^{t}_{j}) /s}(x,y)+\sum_{i\in T^{tj}_{s}\setminus \{\alpha^{tj}_{s}, \beta^{tj}_{s}\}}\mathscr{T}'_{(\widehat{\widehat{P}^{t}_{j}})^{s}_{i}}(x,y)\right)\\
&= &x^{2}\mathscr{T}'_{P\setminus t \setminus s}(x,y)+y^{2}\mathscr{T}'_{P/t\setminus s}(x,y)+xy(\mathscr{T}'_{P\setminus t \setminus s}(x,y)+\mathscr{T}'_{P/t\setminus s}(x,y))\\
& &+x\left(\sum_{j\in T_{t}\setminus \{\alpha_{t}, \beta_{t}\}}\mathscr{T}'_{(\widehat{P}^{t}_{j})\setminus s}(x,y)-\mathscr{T}'_{P\setminus t \setminus s}(x,y)\right)\\
& &+y\left(\sum_{j\in T_{t}\setminus \{\alpha_{t}, \beta_{t}\}}\mathscr{T}'_{(\widehat{P}^{t}_{j})/s}(x,y)-\mathscr{T}'_{P/t\setminus s}(x,y)\right)\\
& &+\sum_{j\in T_{t}\setminus \{\alpha_{t}, \beta_{t}\}}\sum_{i\in T^{tj}_{s}\setminus \{\alpha^{tj}_{s}, \beta^{tj}_{s}\}}\mathscr{T}'_{(\widehat{\widehat{P}^{t}_{j}})^{s}_{i}}(x,y).
\end{eqnarray*}

If first deal with $s$, and then deal with $t$, then by Sub-Case 3.1,
\begin{eqnarray*}
&&	\mathscr{T}'_{P}(x,y)\\
&=&	 x^{2}\mathscr{T}'_{P\setminus s \setminus t}(x,y)+y^{2}\mathscr{T}'_{P/s/t}(x,y)+xy(\mathscr{T}'_{P\setminus s /t}(x,y)+\mathscr{T}'_{P/s\setminus t}(x,y))\\
& &+x\left(\sum_{j\in T^{s\alpha_{s}}_{t}\setminus \{\alpha^{s\alpha_{s}}_{t}, \beta^{s\alpha_{s}}_{t}\}}\mathscr{T}'_{(\widehat{P\setminus s})^{t}_{j}}(x,y)+\sum_{i\in T_{s}\setminus \{\alpha_{s}, \beta_{s}\}}\mathscr{T}'_{(\widehat{P}^{s}_{i})\setminus t}(x,y)\right)\\
& &+y\left(\sum_{j\in T^{s\beta_{s}}_{t}\setminus \{\alpha^{s\beta_{s}}_{t}, \beta^{s\beta_{s}}_{t}\}}\mathscr{T}'_{(\widehat{P/s})^{t}_{j}}(x,y)+\sum_{i\in T_{s}\setminus \{\alpha_{s}, \beta_{s}\}}\mathscr{T}'_{(\widehat{P}^{s}_{i})/t}(x,y)\right)\\
& &+\sum_{i\in T_{s}\setminus \{\alpha_{s}, \beta_{s}\}}\sum_{j\in T^{si}_{t}\setminus \{\alpha^{si}_{t}, \beta^{si}_{t}\}}\mathscr{T}'_{(\widehat{\widehat{P}^{s}_{i}})^{t}_{j}}(x,y).
\end{eqnarray*}
By the analysis for $t$ and $s$ in polymatroids (see Table 2) and induction hypothesis, the conclusion holds.
\begin{table}[htp]
\begin{center}
\renewcommand\arraystretch{1.2}\footnotesize
\caption{A comparison of the polymatroid Tutte polynomial obtained by dealing with $s,t$.}
 \begin{tabular}{|c |c|c|c|c|c|c|c|c|c|c |c|c|c|c|}
  \hline
 \text{first}~$t$  \text{and then}~$s$ &$a_{t}$& $a_{s}$& \text{first}~$s$  \text{and then}~$t$\\
 \hline
 $\textbf{a}\in P\setminus t\setminus s$& $\alpha_{t}$ & $\beta_{s}$ & $\textbf{a}\in P/s\setminus t$\\
 \hline
 $\textbf{a}\in P/t\setminus s$&$\beta_{t}$& $\alpha_{s}$& $\textbf{a}\in P\setminus s /t$\\
 \hline
 $\textbf{a}\in(\widehat{P}^{t}_{j})\setminus s$, \text{where}&$j$& $\alpha_{s}$& $\textbf{a}\in (\widehat{P\setminus s})^{t}_{j}$, \text{where}\\

$j\in T^{s\alpha_{s}}_{t}\setminus \{\alpha^{s\alpha_{s}}_{t}, \beta^{s\alpha_{s}}_{t}\}$&& &$j\in T^{s\alpha_{s}}_{t}\setminus \{\alpha^{s\alpha_{s}}_{t}, \beta^{s\alpha_{s}}_{t}\}$ \\
 \hline
 $\textbf{a}\in(\widehat{P}^{t}_{j})\setminus s$,~\text{where}~$j=$&$f([n]\setminus \{s\})$& $\alpha_{s}$& $\textbf{a}\in P\setminus s\setminus t$\\

$f([n]\setminus\{s\})-f([n]\setminus \{s,t\})$&$-f([n]\setminus \{s,t\})$& & \\
 \hline
 $\textbf{a}\in(\widehat{P}^{t}_{j})\setminus s$, \text{where}~$\alpha_{t}<j<$&$j$& $f([n])-j$& $\textbf{a}\in(\widehat{P}^{s}_{i})\setminus t$, \text{where}~$i=$\\
$f([n]\setminus\{s\})-f([n]\setminus \{s,t\})$& &$-f([n]\setminus \{s,t\})$ & $f([n])-j-f([n]\setminus \{s,t\})$\\
 \hline
 $\textbf{a}\in(\widehat{P}^{t}_{j})/s$, \text{where}&$j$& $\beta_{s}$& $\textbf{a}\in (\widehat{P/s})^{t}_{j}$, \text{where}\\

$j\in T^{s\beta_{s}}_{t}\setminus \{\alpha^{s\beta_{s}}_{t}, \beta^{s\beta_{s}}_{t}\}$&& &$j\in T^{s\beta_{s}}_{t}\setminus \{\alpha^{s\beta_{s}}_{t}, \beta^{s\beta_{s}}_{t}\}$ \\
 \hline
 $\textbf{a}\in(\widehat{P}^{t}_{j})/s$,~\text{where}~&$f(\{s,t\})$& $\beta_{s}$& $\textbf{a}\in P/s/t$\\

$j=f(\{s,t\})-f(\{s\})$&$-f(\{s\})$& & \\
 \hline
 $\textbf{a}\in(\widehat{P}^{t}_{j})/s$, \text{where}&$j$& $f(\{s,t\})-j$& $\textbf{a}\in(\widehat{P}^{s}_{i})/t$, \text{where}\\
$f(\{s,t\})-f(\{s\})<j<\beta_{t}$& && $i=f(\{s,t\})-j$\\
 \hline
\end{tabular}
\end{center}
\end{table}

\textbf{(3-2-3)} Assume that $f([n]\setminus \{t\})+f(\{s,t\})= f([n])+f(\{s\})$ and $f([n]\setminus \{s,t\})+f(\{t\})= f([n]\setminus \{s\})$. It is similar to (3-2-2).

\textbf{(3-2-4)} Assume that $f([n]\setminus \{s\})+f(\{s,t\})= f([n])+f(\{t\})$ and $f([n]\setminus \{t\})+f(\{s,t\})= f([n])+f(\{s\})$.
In this case, $\alpha^{t\beta_{t}}_{s}=\beta^{t\beta_{t}}_{s}$, $\alpha^{t\alpha_{t}}_{s}<\beta^{t\alpha_{t}}_{s}$, $\alpha^{s\alpha_{s}}_{t}<\beta^{s\alpha_{s}}_{t}$ and $\alpha^{s\beta_{s}}_{t}=\beta^{s\beta_{s}}_{t}$. Moreover, $P/ t\setminus s=P/ t/s=P/ s/t=P/ s\setminus t$.

If first deal with $t$, and then deal with $s$, then
\begin{eqnarray*}
&&	\mathscr{T}'_{P}(x,y)\\
&=&	x\mathscr{T}'_{P\setminus t}(x,y)+y\mathscr{T}'_{P/t}(x,y)+\sum_{j\in T_{t}\setminus \{\alpha_{t}, \beta_{t}\}}\mathscr{T}'_{\widehat{P}^{t}_{j}}(x,y)\\
&= &x\left(x\mathscr{T}'_{P\setminus t \setminus s}(x,y)+y\mathscr{T}'_{P\setminus t /s}(x,y)+\sum_{i\in T^{t\alpha_{t}}_{s}\setminus \{\alpha^{t\alpha_{t}}_{s}, \beta^{t\alpha_{t}}_{s}\}}\mathscr{T}'_{(\widehat{P\setminus t})^{s}_{i}}(x,y)\right)\\
& &+y(x+y-1)\mathscr{T}'_{P/t\setminus s}(x,y)\\
& &+\sum_{j\in T_{t}\setminus \{\alpha_{t}, \beta_{t}\}}\left(x\mathscr{T}'_{(\widehat{P}^{t}_{j})\setminus s}(x,y)+y\mathscr{T}'_{(\widehat{P}^{t}_{j}) /s}(x,y)+\sum_{i\in T^{tj}_{s}\setminus \{\alpha^{tj}_{s}, \beta^{tj}_{s}\}}\mathscr{T}'_{(\widehat{\widehat{P}^{t}_{j}})^{s}_{i}}(x,y)\right)\\
&= &x^{2}\mathscr{T}'_{P\setminus t \setminus s}(x,y)+y^{2}\mathscr{T}'_{P/t\setminus s}(x,y)+xy(\mathscr{T}'_{P\setminus t/s}(x,y)+\mathscr{T}'_{P/t\setminus s}(x,y))\\
& &+x\left(\sum_{j\in T_{t}\setminus \{\alpha_{t}, \beta_{t}\}}\mathscr{T}'_{(\widehat{P}^{t}_{j})\setminus s}(x,y)+\sum_{i\in T^{t\alpha_{t}}_{s}\setminus \{\alpha^{t\alpha_{t}}_{s}, \beta^{t\alpha_{t}}_{s}\}}\mathscr{T}'_{(\widehat{P\setminus t})^{s}_{i}}(x,y)\right)\\
& &+y\left(\sum_{j\in T_{t}\setminus \{\alpha_{t}, \beta_{t}\}}\mathscr{T}'_{(\widehat{P}^{t}_{j})/s}(x,y)-\mathscr{T}'_{P/t\setminus s}(x,y)\right)\\
& &+\sum_{j\in T_{t}\setminus \{\alpha_{t}, \beta_{t}\}}\sum_{i\in T^{tj}_{s}\setminus \{\alpha^{tj}_{s}, \beta^{tj}_{s}\}}\mathscr{T}'_{(\widehat{\widehat{P}^{t}_{j}})^{s}_{i}}(x,y).
\end{eqnarray*}

If first deal with $s$, and then deal with $t$, then
\begin{eqnarray*}
&&	\mathscr{T}'_{P}(x,y)\\
&=&	 x^{2}\mathscr{T}'_{P\setminus s \setminus t}(x,y)+y^{2}\mathscr{T}'_{P/s\setminus t}(x,y)+xy(\mathscr{T}'_{P\setminus s /t}(x,y)+\mathscr{T}'_{P/s\setminus t}(x,y))\\
& &+x\left(\sum_{j\in T^{s\alpha_{s}}_{t}\setminus \{\alpha^{s\alpha_{s}}_{t}, \beta^{s\alpha_{s}}_{t}\}}\mathscr{T}'_{(\widehat{P\setminus s})^{t}_{j}}(x,y)+\sum_{i\in T_{s}\setminus \{\alpha_{s}, \beta_{s}\}}\mathscr{T}'_{(\widehat{P}^{s}_{i})\setminus t}(x,y)\right)\\
& &+y\left(\sum_{i\in T_{s}\setminus \{\alpha_{s}, \beta_{s}\}}\mathscr{T}'_{(\widehat{P}^{s}_{i})/t}(x,y)-\mathscr{T}'_{P/s\setminus t}(x,y)\right)\\
& &+\sum_{i\in T_{s}\setminus \{\alpha_{s}, \beta_{s}\}}\sum_{j\in T^{si}_{t}\setminus \{\alpha^{si}_{t}, \beta^{si}_{t}\}}\mathscr{T}'_{(\widehat{\widehat{P}^{s}_{i}})^{t}_{j}}(x,y).
\end{eqnarray*}
By the analysis for $t$ and $s$ in polymatroids (see Table 3) and induction hypothesis, the conclusion holds.
\begin{table}[htp]
\begin{center}
\renewcommand\arraystretch{1.2}\footnotesize
\caption{A comparison of the polymatroid Tutte polynomial obtained by dealing with $s,t$.}
 \begin{tabular}{|c |c|c|c|c|c|c|c|c|c|c |c|c|c|c|}
  \hline
 \text{first}~$t$  \text{and then}~$s$ &$a_{t}$& $a_{s}$& \text{first}~$s$  \text{and then}~$t$\\
 \hline
 $\textbf{a}\in P\setminus t\setminus s$& $\alpha_{t}$ & $f([n]\setminus \{t\})-$ & $\textbf{a}\in(\widehat{P}^{s}_{i})\setminus t$, \text{where}~$i=$\\
 & &$f([n]\setminus \{s,t\})$ &$f([n]\setminus\{t\})-f([n]\setminus \{s,t\})$\\
 \hline
$\textbf{a}\in P\setminus t /s$ & $\alpha_{t}$ & $\beta_{s}$ & $\textbf{a}\in P/s\setminus t$\\
 \hline
 $\textbf{a}\in P/t\setminus s$&$\beta_{t}$& $\alpha_{s}$& $\textbf{a}\in P\setminus s /t$\\
 \hline
 $\textbf{a}\in (\widehat{P\setminus t})^{s}_{i}$, &$\alpha_{t}$& $i$& $\textbf{a}\in(\widehat{P}^{s}_{i})\setminus t$, \text{where}~ \\

\text{where}~$i\in T^{t\alpha_{t}}_{s}\setminus \{\alpha^{t\alpha_{t}}_{s}, \beta^{t\alpha_{t}}_{s}\}$& & &$i\in T^{t\alpha_{t}}_{s}\setminus \{\alpha^{t\alpha_{t}}_{s}, \beta^{t\alpha_{t}}_{s}\}$\\
 \hline
 $\textbf{a}\in(\widehat{P}^{t}_{j})\setminus s$, \text{where}&$j$& $\alpha_{s}$& $\textbf{a}\in (\widehat{P\setminus s})^{t}_{j}$, \text{where}\\

$j\in T^{s\alpha_{s}}_{t}\setminus \{\alpha^{s\alpha_{s}}_{t}, \beta^{s\alpha_{s}}_{t}\}$&& &$j\in T^{s\alpha_{s}}_{t}\setminus \{\alpha^{s\alpha_{s}}_{t}, \beta^{s\alpha_{s}}_{t}\}$ \\
 \hline
 $\textbf{a}\in(\widehat{P}^{t}_{j})\setminus s$,~\text{where}~$j=$&$f([n]\setminus \{s\})$& $\alpha_{s}$& $\textbf{a}\in P\setminus s\setminus t$\\

$f([n]\setminus\{s\})-f([n]\setminus \{s,t\})$&$-f([n]\setminus \{s,t\})$& & \\
 \hline
 $\textbf{a}\in(\widehat{P}^{t}_{j})\setminus s$, \text{where}~$\alpha_{t}<j<$&$j$& $f([n])-j-$& $\textbf{a}\in(\widehat{P}^{s}_{i})\setminus t$, \text{where}~$i=$\\
$f([n]\setminus\{s\})-f([n]\setminus \{s,t\})$& &$f([n]\setminus \{s,t\})$ & $f([n])-j-f([n]\setminus \{s,t\})$\\
 \hline
 $\textbf{a}\in(\widehat{P}^{t}_{j})/s$, \text{where}&$j$& $f(\{s,t\})-j$& $\textbf{a}\in(\widehat{P}^{s}_{i})/t$, \text{where}\\
$f(\{s,t\})-f(\{s\})<j<\beta_{t}$& && $i=f(\{s,t\})-j$\\
 \hline
\end{tabular}
\end{center}
\end{table}

\textbf{(3-2-5)} Assume that $f([n]\setminus \{s,t\})+f(\{s\})= f([n]\setminus \{t\})$ and $f([n]\setminus \{s,t\})+f(\{t\})=f([n]\setminus \{s\})$.
Then $\alpha^{t\alpha_{t}}_{s}=\beta^{t\alpha_{t}}_{s}$, $\alpha^{s\alpha_{s}}_{t}=\beta^{s\alpha_{s}}_{t}$, $\alpha^{t\beta_{t}}_{s}<\beta^{t\beta_{t}}_{s}$ and $\alpha^{s\beta_{s}}_{t}<\beta^{s\beta_{s}}_{t}$. Moreover, $P\setminus t/ s=P\setminus t\setminus s=P\setminus s\setminus t=P\setminus s/ t$.

If first deal with $t$, and then deal with $s$, then
\begin{eqnarray*}
&&	\mathscr{T}'_{P}(x,y)\\
&=&	 x\mathscr{T}'_{P\setminus t}(x,y)+y\mathscr{T}'_{P/t}(x,y)+\sum_{j\in T_{t}\setminus \{\alpha_{t}, \beta_{t}\}}\mathscr{T}'_{\widehat{P}^{t}_{j}}(x,y)\\
	&= &x(x+y-1)\mathscr{T}'_{P\setminus t \setminus s}(x,y)\\
& &+y\left(x\mathscr{T}'_{P/t\setminus s}(x,y)+y\mathscr{T}'_{P/t/s}(x,y)+\sum_{i\in T^{t\beta_{t}}_{s}\setminus \{\alpha^{t\beta_{t}}_{s}, \beta^{t\beta_{t}}_{s}\}}\mathscr{T}'_{(\widehat{P/t})^{s}_{i}}(x,y)\right)\\
& &+\sum_{j\in T_{t}\setminus \{\alpha_{t}, \beta_{t}\}}\left(x\mathscr{T}'_{(\widehat{P}^{t}_{j})\setminus s}(x,y)+y\mathscr{T}'_{(\widehat{P}^{t}_{j}) /s}(x,y)+\sum_{i\in T^{tj}_{s}\setminus \{\alpha^{tj}_{s}, \beta^{tj}_{s}\}}\mathscr{T}'_{(\widehat{\widehat{P}^{t}_{j}})^{s}_{i}}(x,y)\right)\\
&= &x^{2}\mathscr{T}'_{P\setminus t \setminus s}(x,y)+y^{2}\mathscr{T}'_{P/t/s}(x,y)+xy(\mathscr{T}'_{P\setminus t\setminus s}(x,y)+\mathscr{T}'_{P/t\setminus s}(x,y))\\
& &+x\left(\sum_{j\in T_{t}\setminus \{\alpha_{t}, \beta_{t}\}}\mathscr{T}'_{(\widehat{P}^{t}_{j})\setminus s}(x,y)-\mathscr{T}'_{P\setminus t \setminus s}(x,y)\right)\\
& &+y\left(\sum_{i\in T^{t\beta_{t}}_{s}\setminus \{\alpha^{t\beta_{t}}_{s}, \beta^{t\beta_{t}}_{s}\}}\mathscr{T}'_{(\widehat{P/t})^{s}_{i}}(x,y)+\sum_{j\in T_{t}\setminus \{\alpha_{t}, \beta_{t}\}}\mathscr{T}'_{(\widehat{P}^{t}_{j})/s}(x,y)\right)\\
& &+\sum_{j\in T_{t}\setminus \{\alpha_{t}, \beta_{t}\}}\sum_{i\in T^{tj}_{s}\setminus \{\alpha^{tj}_{s}, \beta^{tj}_{s}\}}\mathscr{T}'_{(\widehat{\widehat{P}^{t}_{j}})^{s}_{i}}(x,y).
\end{eqnarray*}
If first deal with $s$, and then deal with $t$, then
\begin{eqnarray*}
&&	\mathscr{T}'_{P}(x,y)\\
&=&	 x^{2}\mathscr{T}'_{P\setminus s \setminus t}(x,y)+y^{2}\mathscr{T}'_{P/s/t}(x,y)+xy(\mathscr{T}'_{P\setminus s \setminus t}(x,y)+\mathscr{T}'_{P/s\setminus t}(x,y))\\
& &+x\left(\sum_{i\in T_{s}\setminus \{\alpha_{s}, \beta_{s}\}}\mathscr{T}'_{(\widehat{P}^{s}_{i})\setminus t}(x,y)-\mathscr{T}'_{P\setminus s \setminus t}(x,y)\right)\\
& &+y\left(\sum_{j\in T^{s\beta_{s}}_{t}\setminus \{\alpha^{s\beta_{s}}_{t}, \beta^{s\beta_{s}}_{t}\}}\mathscr{T}'_{(\widehat{P/s})^{t}_{j}}(x,y)+\sum_{i\in T_{s}\setminus \{\alpha_{s}, \beta_{s}\}}\mathscr{T}'_{(\widehat{P}^{s}_{i})/t}(x,y)\right)\\
& &+\sum_{i\in T_{s}\setminus \{\alpha_{s}, \beta_{s}\}}\sum_{j\in T^{si}_{t}\setminus \{\alpha^{si}_{t}, \beta^{si}_{t}\}}\mathscr{T}'_{(\widehat{\widehat{P}^{s}_{i}})^{t}_{j}}(x,y).
\end{eqnarray*}

By the analysis for $t$ and $s$ in polymatroids (see Table 4) and induction hypothesis, the conclusion holds.
\begin{table}[htp]
\begin{center}
\renewcommand\arraystretch{1.2} \footnotesize
\caption{A comparison of the polymatroid Tutte polynomial obtained by dealing with $s,t$.}
 \begin{tabular}{|c |c|c|c|c|c|c|c|c|c|c |c|c|c|c|}
  \hline
 \text{first}~$t$  \text{and then}~$s$ &$a_{t}$& $a_{s}$& \text{first}~$s$  \text{and then}~$t$\\
 \hline
 $\textbf{a}\in P\setminus t\setminus s$& $\alpha_{t}$ & $\beta_{s}$ & $\textbf{a}\in P/s\setminus t$\\
 \hline
 $\textbf{a}\in P/t/s$&$\beta_{t}$& $f(\{s,t\})-$& $\textbf{a}\in(\widehat{P}^{t}_{j})\setminus s$, \text{where}\\
 & &$f(\{t\})$ &$i=f(\{s,t\})$\\
 \hline
 $\textbf{a}\in P/t\setminus s$&$\beta_{t}$& $\alpha_{s}$& $\textbf{a}\in P\setminus s\setminus t$\\

 \hline
 $\textbf{a}\in(\widehat{P}^{t}_{j})\setminus s$, \text{where}~$\alpha_{t}<j<$&$j$& $f([n])-j-$& $\textbf{a}\in(\widehat{P}^{s}_{i})\setminus t$, \text{where}~$i=$\\
$f([n]\setminus\{s\})-f([n]\setminus \{s,t\})$& &$f([n]\setminus \{s,t\})$ & $f([n])-j-f([n]\setminus \{s,t\})$\\
 \hline
 $\textbf{a}\in (\widehat{P/t})^{s}_{i}$, &$\beta_{t}$& $i$& $\textbf{a}\in(\widehat{P}^{s}_{i})/t$, \text{where}~ \\

\text{where}~$i\in T^{t\beta_{t}}_{s}\setminus \{\alpha^{t\beta_{t}}_{s}, \beta^{t\beta_{t}}_{s}\}$& & &$i\in T^{t\beta_{t}}_{s}\setminus \{\alpha^{t\beta_{t}}_{s}, \beta^{t\beta_{t}}_{s}\}$\\
 \hline
 $\textbf{a}\in(\widehat{P}^{t}_{j})/s$, \text{where}&$j$& $\beta_{s}$& $\textbf{a}\in (\widehat{P/s})^{t}_{j}$, \text{where}\\

$j\in T^{s\beta_{s}}_{t}\setminus \{\alpha^{s\beta_{s}}_{t}, \beta^{s\beta_{s}}_{t}\}$&& &$j\in T^{s\beta_{s}}_{t}\setminus \{\alpha^{s\beta_{s}}_{t}, \beta^{s\beta_{s}}_{t}\}$ \\
 \hline
 $\textbf{a}\in(\widehat{P}^{t}_{j})/s$,~\text{where}~&$f(\{s,t\})$& $\beta_{s}$& $\textbf{a}\in P/s/t$\\

$j=f(\{s,t\})-f(\{s\})$&$-f(\{s\})$& & \\
 \hline
 $\textbf{a}\in(\widehat{P}^{t}_{j})/s$, \text{where}&$j$& $f(\{s,t\})-j$& $\textbf{a}\in(\widehat{P}^{s}_{i})/t$, \text{where}\\
$f(\{s,t\})-f(\{s\})<j<\beta_{t}$& && $i=f(\{s,t\})-j$\\
 \hline
\end{tabular}
\end{center}
\end{table}

\textbf{(3-2-6)} Assume that $f([n]\setminus \{s\})+f(\{s,t\})= f([n])+f(\{t\})$, $f([n]\setminus \{s,t\})+f(\{s\})> f([n]\setminus \{t\})$ and $f([n]\setminus \{t\})+f(\{s,t\})>f([n])+f(\{s\})$.
In this case, $\alpha^{t\beta_{t}}_{s}=\beta^{t\beta_{t}}_{s}$, $\alpha^{t\alpha_{t}}_{s}<\beta^{t\alpha_{t}}_{s}$, $\alpha^{s\alpha_{s}}_{t}<\beta^{s\alpha_{s}}_{t}$ and $\alpha^{s\beta_{s}}_{t}<\beta^{s\beta_{s}}_{t}$.

If first deal with $t$, and then deal with $s$, then
\begin{eqnarray*}
&&	\mathscr{T}'_{P}(x,y)\\
&=&	x^{2}\mathscr{T}'_{P\setminus t \setminus s}(x,y)+y^{2}\mathscr{T}'_{P/t\setminus s}(x,y)+xy(\mathscr{T}'_{P\setminus t/s}(x,y)+\mathscr{T}'_{P/t\setminus s}(x,y))\\
& &+x\left(\sum_{j\in T_{t}\setminus \{\alpha_{t}, \beta_{t}\}}\mathscr{T}'_{(\widehat{P}^{t}_{j})\setminus s}(x,y)+\sum_{i\in T^{t\alpha_{t}}_{s}\setminus \{\alpha^{t\alpha_{t}}_{s}, \beta^{t\alpha_{t}}_{s}\}}\mathscr{T}'_{(\widehat{P\setminus t})^{s}_{i}}(x,y)\right)\\
& &+y\left(\sum_{j\in T_{t}\setminus \{\alpha_{t}, \beta_{t}\}}\mathscr{T}'_{(\widehat{P}^{t}_{j})/s}(x,y)-\mathscr{T}'_{P/t\setminus s}(x,y)\right)\\
& &+\sum_{j\in T_{t}\setminus \{\alpha_{t}, \beta_{t}\}}\sum_{i\in T^{tj}_{s}\setminus \{\alpha^{tj}_{s}, \beta^{tj}_{s}\}}\mathscr{T}'_{(\widehat{\widehat{P}^{t}_{j}})^{s}_{i}}(x,y).
\end{eqnarray*}

If first deal with $s$, and then deal with $t$, then by Sub-Case 3.1,
\begin{eqnarray*}
&&	\mathscr{T}'_{P}(x,y)\\
&=&	 x^{2}\mathscr{T}'_{P\setminus s \setminus t}(x,y)+y^{2}\mathscr{T}'_{P/s/t}(x,y)+xy(\mathscr{T}'_{P\setminus s /t}(x,y)+\mathscr{T}'_{P/s\setminus t}(x,y))\\
& &+x\left(\sum_{j\in T^{s\alpha_{s}}_{t}\setminus \{\alpha^{s\alpha_{s}}_{t}, \beta^{s\alpha_{s}}_{t}\}}\mathscr{T}'_{(\widehat{P\setminus s})^{t}_{j}}(x,y)+\sum_{i\in T_{s}\setminus \{\alpha_{s}, \beta_{s}\}}\mathscr{T}'_{(\widehat{P}^{s}_{i})\setminus t}(x,y)\right)\\
& &+y\left(\sum_{j\in T^{s\beta_{s}}_{t}\setminus \{\alpha^{s\beta_{s}}_{t}, \beta^{s\beta_{s}}_{t}\}}\mathscr{T}'_{(\widehat{P/s})^{t}_{j}}(x,y)+\sum_{i\in T_{s}\setminus \{\alpha_{s}, \beta_{s}\}}\mathscr{T}'_{(\widehat{P}^{s}_{i})/t}(x,y)\right)\\
& &+\sum_{i\in T_{s}\setminus \{\alpha_{s}, \beta_{s}\}}\sum_{j\in T^{si}_{t}\setminus \{\alpha^{si}_{t}, \beta^{si}_{t}\}}\mathscr{T}'_{(\widehat{\widehat{P}^{s}_{i}})^{t}_{j}}(x,y).
\end{eqnarray*}

Note that $P/t\setminus s=P/t/s$. By the analysis for $t$ and $s$ in polymatroids (see Table 5) and induction hypothesis, the conclusion holds.
\begin{table}[htp]
\begin{center}
\renewcommand\arraystretch{1.2}\footnotesize
\caption{A comparison of the polymatroid Tutte polynomial obtained by dealing with $s,t$.}
 \begin{tabular}{|c |c|c|c|c|c|c|c|c|c|c |c|c|c|c|}
  \hline
 \text{first}~$t$  \text{and then}~$s$ &$a_{t}$& $a_{s}$& \text{first}~$s$  \text{and then}~$t$\\
 \hline
 $\textbf{a}\in P\setminus t\setminus s$& $\alpha_{t}$ & $f([n]\setminus \{t\})-$ & $\textbf{a}\in(\widehat{P}^{s}_{i})\setminus t$, \text{where}~$i=$\\
 & &$f([n]\setminus \{s,t\})$ &$f([n]\setminus\{t\})-f([n]\setminus \{s,t\})$\\
 \hline
$\textbf{a}\in P\setminus t /s$ & $\alpha_{t}$ & $\beta_{s}$ & $\textbf{a}\in P/s\setminus t$\\
 \hline
 $\textbf{a}\in P/t\setminus s$&$\beta_{t}$& $\alpha_{s}$& $\textbf{a}\in P\setminus s /t$\\
 \hline
 $\textbf{a}\in (\widehat{P\setminus t})^{s}_{i}$, &$\alpha_{t}$& $i$& $\textbf{a}\in(\widehat{P}^{s}_{i})\setminus t$, \text{where}~ \\

\text{where}~$i\in T^{t\alpha_{t}}_{s}\setminus \{\alpha^{t\alpha_{t}}_{s}, \beta^{t\alpha_{t}}_{s}\}$& & &$i\in T^{t\alpha_{t}}_{s}\setminus \{\alpha^{t\alpha_{t}}_{s}, \beta^{t\alpha_{t}}_{s}\}$\\
 \hline
 $\textbf{a}\in(\widehat{P}^{t}_{j})\setminus s$, \text{where}&$j$& $\alpha_{s}$& $\textbf{a}\in (\widehat{P\setminus s})^{t}_{j}$, \text{where}\\

$j\in T^{s\alpha_{s}}_{t}\setminus \{\alpha^{s\alpha_{s}}_{t}, \beta^{s\alpha_{s}}_{t}\}$&& &$j\in T^{s\alpha_{s}}_{t}\setminus \{\alpha^{s\alpha_{s}}_{t}, \beta^{s\alpha_{s}}_{t}\}$ \\
 \hline
 $\textbf{a}\in(\widehat{P}^{t}_{j})\setminus s$,~\text{where}~$j=$&$f([n]\setminus \{s\})$& $\alpha_{s}$& $\textbf{a}\in P\setminus s\setminus t$\\

$f([n]\setminus\{s\})-f([n]\setminus \{s,t\})$&$-f([n]\setminus \{s,t\})$& & \\
 \hline
 $\textbf{a}\in(\widehat{P}^{t}_{j})\setminus s$, \text{where}~$\alpha_{t}<j<$&$j$& $f([n])-j-$& $\textbf{a}\in(\widehat{P}^{s}_{i})\setminus t$, \text{where}~$i=$\\
$f([n]\setminus\{s\})-f([n]\setminus \{s,t\})$& &$f([n]\setminus \{s,t\})$ & $f([n])-j-f([n]\setminus \{s,t\})$\\
 \hline
 $\textbf{a}\in(\widehat{P}^{t}_{j})/s$, \text{where}&$j$& $\beta_{s}$& $\textbf{a}\in (\widehat{P/s})^{t}_{j}$, \text{where}\\

$j\in T^{s\beta_{s}}_{t}\setminus \{\alpha^{s\beta_{s}}_{t}, \beta^{s\beta_{s}}_{t}\}$&& &$j\in T^{s\beta_{s}}_{t}\setminus \{\alpha^{s\beta_{s}}_{t}, \beta^{s\beta_{s}}_{t}\}$ \\
 \hline
 $\textbf{a}\in(\widehat{P}^{t}_{j})/s$,~\text{where}~&$f(\{s,t\})$& $\beta_{s}$& $\textbf{a}\in P/s/t$\\

$j=f(\{s,t\})-f(\{s\})$&$-f(\{s\})$& & \\
 \hline
 $\textbf{a}\in(\widehat{P}^{t}_{j})/s$, \text{where}&$j$& $f(\{s,t\})-j$& $\textbf{a}\in(\widehat{P}^{s}_{i})/t$, \text{where}\\
$f(\{s,t\})-f(\{s\})<j<\beta_{t}$& && $i=f(\{s,t\})-j$\\
 \hline
\end{tabular}
\end{center}
\end{table}

\textbf{(3-2-7)} Assume that $f([n]\setminus \{t\})+f(\{s,t\})= f([n])+f(\{s\})$, $f([n]\setminus \{s,t\})+f(\{t\})> f([n]\setminus \{s\})$ and $f([n]\setminus \{s\})+f(\{s,t\})>f([n])+f(\{t\})$. It is similar to (3-2-6).

\textbf{(3-2-8)} Assume that $f([n]\setminus \{s,t\})+f(\{s\})= f([n]\setminus \{t\})$,  $f([n]\setminus \{s\})+f(\{s,t\})> f([n])+f(\{t\})$ and $f([n]\setminus \{s,t\})+f(\{t\})>f([n]\setminus \{s\})$.
Then $\alpha^{t\alpha_{t}}_{s}=\beta^{t\alpha_{t}}_{s}$, $\alpha^{s\alpha_{s}}_{t}<\beta^{s\alpha_{s}}_{t}$, $\alpha^{t\beta_{t}}_{s}<\beta^{t\beta_{t}}_{s}$ and $\alpha^{s\beta_{s}}_{t}<\beta^{s\beta_{s}}_{t}$. Moreover, $P\setminus t\setminus s=P\setminus t/s$.

If first deal with $t$, and then deal with $s$, then
\begin{eqnarray*}
&&	\mathscr{T}'_{P}(x,y)\\
&=&	 x^{2}\mathscr{T}'_{P\setminus t \setminus s}(x,y)+y^{2}\mathscr{T}'_{P/t/s}(x,y)+xy(\mathscr{T}'_{P\setminus t\setminus s}(x,y)+\mathscr{T}'_{P/t\setminus s}(x,y))\\
& &+x\left(\sum_{j\in T_{t}\setminus \{\alpha_{t}, \beta_{t}\}}\mathscr{T}'_{(\widehat{P}^{t}_{j})\setminus s}(x,y)-\mathscr{T}'_{P\setminus t \setminus s}(x,y)\right)\\
& &+y\left(\sum_{i\in T^{t\beta_{t}}_{s}\setminus \{\alpha^{t\beta_{t}}_{s}, \beta^{t\beta_{t}}_{s}\}}\mathscr{T}'_{(\widehat{P/t})^{s}_{i}}(x,y)+\sum_{j\in T_{t}\setminus \{\alpha_{t}, \beta_{t}\}}\mathscr{T}'_{(\widehat{P}^{t}_{j})/s}(x,y)\right)\\
& &+\sum_{j\in T_{t}\setminus \{\alpha_{t}, \beta_{t}\}}\sum_{i\in T^{tj}_{s}\setminus \{\alpha^{tj}_{s}, \beta^{tj}_{s}\}}\mathscr{T}'_{(\widehat{\widehat{P}^{t}_{j}})^{s}_{i}}(x,y).
\end{eqnarray*}

If first deal with $s$, and then deal with $t$, then by Sub-Case 3.1,
\begin{eqnarray*}
&&	\mathscr{T}'_{P}(x,y)\\
&=&	 x^{2}\mathscr{T}'_{P\setminus s \setminus t}(x,y)+y^{2}\mathscr{T}'_{P/s/t}(x,y)+xy(\mathscr{T}'_{P\setminus s /t}(x,y)+\mathscr{T}'_{P/s\setminus t}(x,y))\\
& &+x\left(\sum_{j\in T^{s\alpha_{s}}_{t}\setminus \{\alpha^{s\alpha_{s}}_{t}, \beta^{s\alpha_{s}}_{t}\}}\mathscr{T}'_{(\widehat{P\setminus s})^{t}_{j}}(x,y)+\sum_{i\in T_{s}\setminus \{\alpha_{s}, \beta_{s}\}}\mathscr{T}'_{(\widehat{P}^{s}_{i})\setminus t}(x,y)\right)\\
& &+y\left(\sum_{j\in T^{s\beta_{s}}_{t}\setminus \{\alpha^{s\beta_{s}}_{t}, \beta^{s\beta_{s}}_{t}\}}\mathscr{T}'_{(\widehat{P/s})^{t}_{j}}(x,y)+\sum_{i\in T_{s}\setminus \{\alpha_{s}, \beta_{s}\}}\mathscr{T}'_{(\widehat{P}^{s}_{i})/t}(x,y)\right)\\
& &+\sum_{i\in T_{s}\setminus \{\alpha_{s}, \beta_{s}\}}\sum_{j\in T^{si}_{t}\setminus \{\alpha^{si}_{t}, \beta^{si}_{t}\}}\mathscr{T}'_{(\widehat{\widehat{P}^{s}_{i}})^{t}_{j}}(x,y).
\end{eqnarray*}
By the analysis for $t$ and $s$ in polymatroids (see Table 6) and induction hypothesis, the conclusion holds.
\begin{table}[htp]
\begin{center}
\renewcommand\arraystretch{1.2}\footnotesize
\caption{A comparison of the polymatroid Tutte polynomial obtained by dealing with $s,t$.}
 \begin{tabular}{|c |c|c|c|c|c|c|c|c|c|c |c|c|c|c|}
  \hline
 \text{first}~$t$  \text{and then}~$s$ &$a_{t}$& $a_{s}$& \text{first}~$s$  \text{and then}~$t$\\
 \hline
 $\textbf{a}\in P\setminus t\setminus s$& $\alpha_{t}$ & $\beta_{s}$ & $\textbf{a}\in P/s\setminus t$\\
 \hline
 $\textbf{a}\in P/t/s$&$\beta_{t}$& $f(\{s,t\})-$& $\textbf{a}\in(\widehat{P}^{t}_{j})\setminus s$, \text{where}\\
 & &$f(\{t\})$ &$i=f(\{s,t\})$\\
 \hline
 $\textbf{a}\in P/t\setminus s$&$\beta_{t}$& $\alpha_{s}$& $\textbf{a}\in P\setminus s /t$\\
 \hline
 $\textbf{a}\in(\widehat{P}^{t}_{j})\setminus s$, \text{where}&$j$& $\alpha_{s}$& $\textbf{a}\in (\widehat{P\setminus s})^{t}_{j}$, \text{where}\\

$j\in T^{s\alpha_{s}}_{t}\setminus \{\alpha^{s\alpha_{s}}_{t}, \beta^{s\alpha_{s}}_{t}\}$&& &$j\in T^{s\alpha_{s}}_{t}\setminus \{\alpha^{s\alpha_{s}}_{t}, \beta^{s\alpha_{s}}_{t}\}$ \\
 \hline
 $\textbf{a}\in(\widehat{P}^{t}_{j})\setminus s$,~\text{where}~$j=$&$f([n]\setminus \{s\})$& $\alpha_{s}$& $\textbf{a}\in P\setminus s\setminus t$\\

$f([n]\setminus\{s\})-f([n]\setminus \{s,t\})$&$-f([n]\setminus \{s,t\})$& & \\
 \hline
 $\textbf{a}\in(\widehat{P}^{t}_{j})\setminus s$, \text{where}~$\alpha_{t}<j<$&$j$& $f([n])-j-$& $\textbf{a}\in(\widehat{P}^{s}_{i})\setminus t$, \text{where}~$i=$\\
$f([n]\setminus\{s\})-f([n]\setminus \{s,t\})$& &$f([n]\setminus \{s,t\})$ & $f([n])-j-f([n]\setminus \{s,t\})$\\
 \hline
 $\textbf{a}\in (\widehat{P/t})^{s}_{i}$, &$\beta_{t}$& $i$& $\textbf{a}\in(\widehat{P}^{s}_{i})/t$, \text{where}~ \\

\text{where}~$i\in T^{t\beta_{t}}_{s}\setminus \{\alpha^{t\beta_{t}}_{s}, \beta^{t\beta_{t}}_{s}\}$& & &$i\in T^{t\beta_{t}}_{s}\setminus \{\alpha^{t\beta_{t}}_{s}, \beta^{t\beta_{t}}_{s}\}$\\
 \hline
 $\textbf{a}\in(\widehat{P}^{t}_{j})/s$, \text{where}&$j$& $\beta_{s}$& $\textbf{a}\in (\widehat{P/s})^{t}_{j}$, \text{where}\\

$j\in T^{s\beta_{s}}_{t}\setminus \{\alpha^{s\beta_{s}}_{t}, \beta^{s\beta_{s}}_{t}\}$&& &$j\in T^{s\beta_{s}}_{t}\setminus \{\alpha^{s\beta_{s}}_{t}, \beta^{s\beta_{s}}_{t}\}$ \\
 \hline
 $\textbf{a}\in(\widehat{P}^{t}_{j})/s$,~\text{where}~&$f(\{s,t\})$& $\beta_{s}$& $\textbf{a}\in P/s/t$\\

$j=f(\{s,t\})-f(\{s\})$&$-f(\{s\})$& & \\
 \hline
 $\textbf{a}\in(\widehat{P}^{t}_{j})/s$, \text{where}&$j$& $f(\{s,t\})-j$& $\textbf{a}\in(\widehat{P}^{s}_{i})/t$, \text{where}\\
$f(\{s,t\})-f(\{s\})<j<\beta_{t}$& && $i=f(\{s,t\})-j$\\
 \hline
\end{tabular}
\end{center}
\end{table}

\textbf{(3-2-9)} Assume that $f([n]\setminus \{s,t\})+f(\{t\})= f([n]\setminus \{s\})$,  $f([n]\setminus \{t\})+f(\{s,t\})> f([n])+f(\{s\})$ and $f([n]\setminus \{s,t\})+f(\{s\})>f([n]\setminus \{t\})$. It is similar to (3-2-8).

Hence, the first claim holds.

Note that both  $\mathscr{T}_{P}(x,y)$ and $\mathscr{T}'_{P}(x,y)$ depend only on $P$, $\mathscr{T}_{P}(x,y)=\mathscr{T}'_{P}(x,y)$ when $n=0$, and $\mathscr{T}_{P}(x,y)$ satisfies the recursive relation of $\mathscr{T}'_{P}(x,y)$. Hence, $\mathscr{T}_{P}(x,y)=\mathscr{T}'_{P}(x,y)$. This completes the proof.
\section*{Acknowledgements}
\noindent

This work was supported by National Natural Science Foundation of China (Nos.~12401462, 12571379 and 12171402), the Natural Science Foundation of Shanxi Province (No.~202403021222034), and Shanxi Key Laboratory of Digital Design and Manufacturing.

\section*{References}
\bibliographystyle{model1b-num-names}
\bibliography{<your-bib-database>}

\begin{thebibliography}{9}
\bibitem{Bernardi}  O. Bernardi, T. K\'{a}lm\'{a}n and A. Postnikov,  Universal Tutte polynomial, Adv. Math. 402 (2022) 108355.

\bibitem{Cameron}   A. Cameron and A. Fink, The Tutte polynomial via lattice point counting, J. Combin. Theory Ser. A 188 (2022) 105584.

\bibitem{Crapo}    H.  Crapo, The Tutte polynomial, Aequationes Math. 3 (1969) 211-229.

\bibitem{Guan4} X. Guan, X. Jin and T. K\'{a}lm\'{a}n, A deletion-contraction formula and monotonicity
properties for the polymatroid Tutte polynomial, Int. Math. Res. Not. 19 (2025) rnaf302.

\bibitem{Guan6} X. Guan and X. Jin, A direct proof of well-definedness for the polymatroid Tutte polynomial, Adv. in Appl. Math. 163 (2025) 102809.

\bibitem{Tutte} W. Tutte, A contribution to the theory of chromatic polynomials, Canad. J. Math. 6 (1954) 80-91.

\bibitem{Whitney}    H. Whitney, On the abstract properties of linear dependence, Amer. J. Math. 57 (1935) 509-533.
\end{thebibliography}

\end{document}